\documentclass[12pt]{article}
\usepackage{microtype} \DisableLigatures{encoding = *, family = * }
\usepackage{subeqn}
\usepackage{graphicx}
\usepackage{ifpdf}
\ifpdf  \fi
\usepackage{amsfonts,amssymb,slashbox}
\usepackage{mathptmx,helvet,courier,makeidx,multicol,footmisc}
\usepackage[numbers]{natbib}
\bibpunct{[}{]}{;}{n}{,}{,}

\newcommand{\commentout}[1]{}

\newcommand{\ba}{\begin{array}}
        \newcommand{\ea}{\end{array}}
\newcommand{\bc}{\begin{center}}
        \newcommand{\ec}{\end{center}}
\newcommand{\bdm}{\begin{displaymath}}
        \newcommand{\edm}{\end{displaymath}}
\newcommand{\bds} {\begin{description}}
        \newcommand{\eds} {\end{description}}
\newcommand{\ben}{\begin{enumerate}}
        \newcommand{\een}{\end{enumerate}}
\newcommand{\beq}{\begin{equation}}
        \newcommand{\eeq}{\end{equation}}
\newcommand{\bfg} {\begin{figure}[h]}
        \newcommand{\efg} {\end{figure}}
\newcommand{\bi} {\begin {itemize}}
        \newcommand{\ei} {\end {itemize}}
\newcommand{\bqn}{\begin{eqnarray}}
        \newcommand{\eqn}{\end{eqnarray}}
\newcommand{\bqs}{\begin{eqnarray*}}
        \newcommand{\eqs}{\end{eqnarray*}}
\newcommand{\bsl} {\begin{slide}[8.8in,6.7in]}
        \newcommand{\esl} {\end{slide}}
\newcommand{\bsq}{\begin{subequations}}
        \newcommand{\esq}{\end{subequations}}       
\newcommand{\bss} {\begin{slide*}[9.3in,6.7in]}
        \newcommand{\ess} {\end{slide*}}
\newcommand{\btb} {\begin {table}}
        \newcommand{\etb} {\end {table}}
\newcommand{\bvb}{\begin{verbatim}}
        \newcommand{\evb}{\end{verbatim}}

\newcommand{\m}{\mbox}

\newcommand {\pd}[2] {{\frac {\partial {#1}} {\partial {#2}}}}

\newcommand{\mat}[1]{{{\left[ \ba #1 \ea \right]}}}
\newcommand{\cas}[1]{{{\left \{ \ba #1 \ea \right. }}}

\newcommand{\reff}[1] {{{Figure \ref {#1}}}}
\newcommand{\refe}[1] {{(\ref {#1})}}

\def\a          {{\alpha}}

\def\pmb#1{\setbox0=\hbox{$#1$}%
   \kern-.025em\copy0\kern-\wd0
   \kern.05em\copy0\kern-\wd0
   \kern-.025em\raise.0433em\box0 }

\def\bfdelta{\pmb \delta}
\def\bfsigma{\pmb \sigma}

\def\bfmu{\pmb \mu}
\def\bfnu{\pmb \nu}

\def\eop{{\hfill $\blacksquare$}}
\def\r{{\rho}}
\newtheorem{theorem}{Theorem}[section]

\newtheorem{lemma}[theorem]{Lemma}
\newtheorem{corollary}[theorem]{Corollary}

\def\b  {{\beta}}

\def\B {{\mathbb{B}}}
\def\F {{\mathbb{F}}}

\def\RS {{\mathbb{RS}}}
\def\FF {{\mathbb{FF}}}

\begin{document}
\title{A Riemann solver for a system of hyperbolic conservation laws at a general road junction} 
\author{Wen-Long Jin \footnote{Department of Civil and Environmental Engineering, California Institute for Telecommunications and Information Technology, Institute of Transportation Studies, 4000 Anteater Instruction and Research Bldg, University of California, Irvine, CA 92697-3600. Tel: 949-824-1672. Fax: 949-824-8385. Email: wjin@uci.edu. Corresponding author}}
\maketitle


\begin{abstract}
The kinematic wave model of traffic flow on a road network is a system of hyperbolic conservation laws, for which the Riemann solver is of physical, analytical, and numerical importance.
In this paper, we present a Riemann solver at a general network junction. In the Riemann solver, we replace the entropy condition in \citep{holden1995unidirection} by a local, discrete flux function used in Cell Transmission Model \citep{daganzo1995ctm}. 
To enable such an entropy condition, which is consistent with fair merging and first-in-first-out diverging rules, 
we enlarge the weak solution space by introducing interior states on a set of measure zero, associated with stationary discontinuities at the junction.
In the demand-supply space, we demonstrate that the Riemann problem is uniquely solved, in the sense that stationary states and, therefore, kinematic waves on all links can be uniquely determined from feasible conditions on both stationary and interior states as well as the entropy condition that prescribes boundary fluxes from interior states.
In addition, the resulting global flux function is the same as the local one. Thus the flux function is both invariant and Godunov. 
 
\end{abstract}

\section{Introduction}

\begin{figure}
\bc
\includegraphics[width=4in]{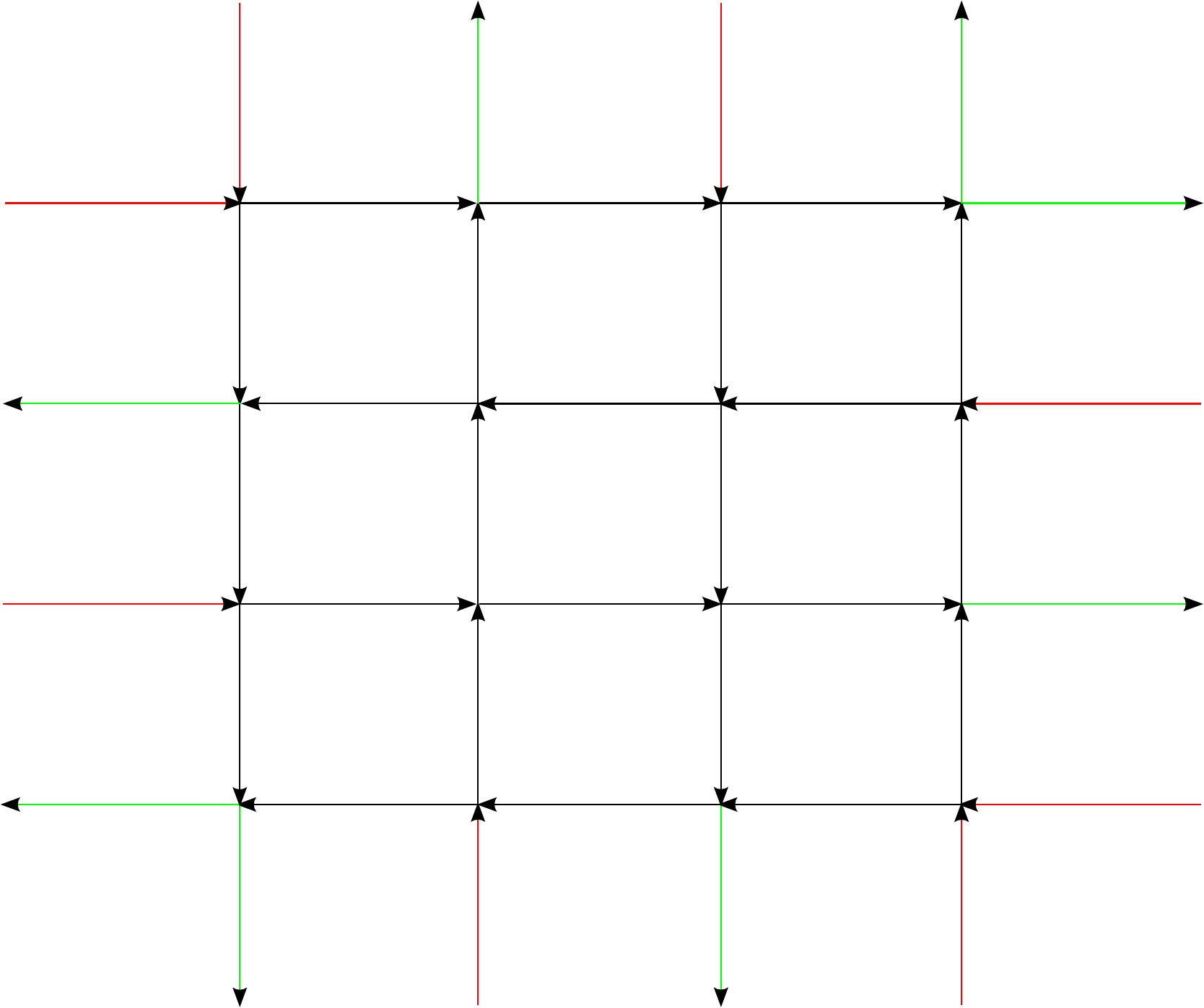}\caption{A grid network}\label{gridnetwork_general}
\ec
\end{figure}

A better understanding of traffic dynamics on a road network is critical for improving safety, mobility, and environmental impacts of modern surface transportation systems \cite{schrank2009mobility}: practically, it is helpful for efficient implementations of ramp metering \cite{papageorgiou2002freeway}, evacuation \cite{Sheffi1982evacuation}, signal control, and other management and control strategies; theoretically, it can yield better network loading models for many other studies \cite{wu1998dnl}.
In a road network, e.g., a grid network shown in \reff{gridnetwork_general}, vehicular traffic dynamics can be described by cellular automata models \cite{nagel1992ca,daganzo2006ca} and car-following models \cite{gazis1961follow} of individual vehicles' movements, fluid dynamic models of continuous car-following behaviors \cite{payne1971PW,whitham1974PW,aw2000arz,zhang2002arz}, the Lighthill-Whitham-Richards (LWR) kinematic wave model \cite{lighthill1955lwr,richards1956lwr}, or regional continuum models \cite{Beckmann1952Transportation,ho2006continuum}.
The traditional LWR model describes traffic dynamics of homogeneous vehicles on a virtually one-lane road as combinations of shock and rarefaction waves and can be analyzed with theories of hyperbolic conservation laws \cite{lax1972shock}.
With the right balance between physical reality and mathematical tractability, kinematic wave models have been successfully extended to study more complicated traffic dynamics of heterogeneous vehicles \cite{benzoni2003populations}, on multi-lane roads \cite{daganzo1997special}, and through network junctions \cite{holden1995unidirection}.

\begin{figure}
\bc\includegraphics[width=3.5in]{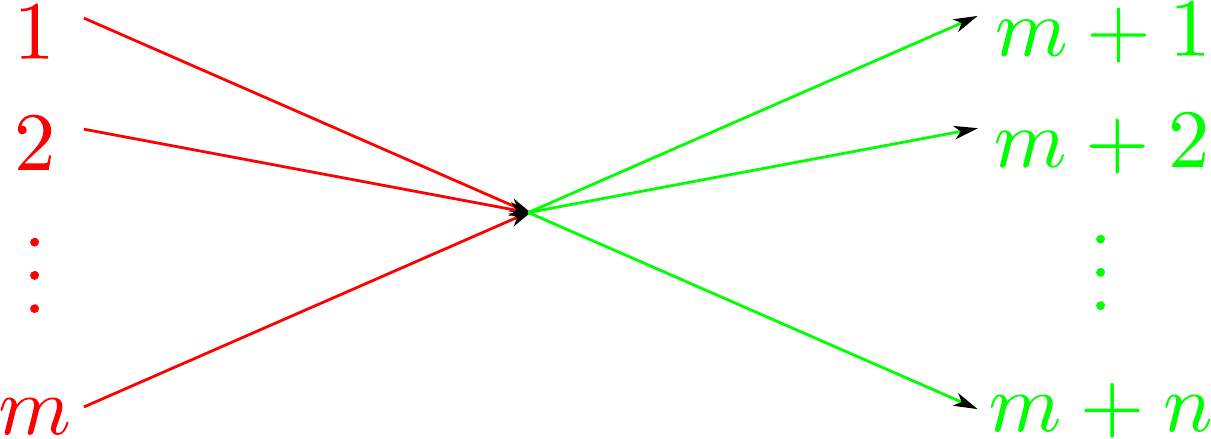}\ec
\caption{A general network junction}\label{general_junction}
\end{figure}

In a road network, such bottlenecks as merging, diverging, and general junctions play a critical role in initiating, propagating, and dissipating traffic congestion. 
Some interesting traffic dynamics can be caused by interactions among these network bottlenecks: for examples, a beltway network can be totally gridlocked  \cite{daganzo1996gridlock}, and periodic oscillations can occur in a diverge-merge network \cite{jin2009network}.
Thus efforts are warranted to develop both physically realistic and mathematically tractable kinematic wave models of  network traffic dynamics.
Since traffic dynamics inside a link can be described by the LWR model and are well understood, the most important component of network kinematic wave models is related to how merging and diverging behaviors would impact the formation of shock and rarefaction waves at a general network junction shown in \reff{general_junction}, which has $m$ upstream links and $n$ downstream links.
In the literature, there have been three lines of research into traffic dynamics through a general network junction. In the line of discrete Cell Transmission Model (CTM) \cite{daganzo1995ctm,lebacque1996godunov}, boundary fluxes through a junction during a time interval are prescribed from adjacent cells' conditions based on macroscopic merging and diverging rules. In the line of continuous models \cite{holden1995unidirection}, shock and rarefaction waves on all links are analytically solved with jump initial conditions by decoupling the Riemann problem at the junction into $m+n$ Riemann problems on individual links. 
In the third line of continuous models \cite{jin2012_network}, the discrete flux functions in CTM are used as decoupling conditions, and it was shown that the Riemann problem can be uniquely solved.
The first two approaches bear their respective limitations: the CTM approach has been verified by empirical observations but cannot be applied to obtain such analytical insights as shock and rarefaction waves at a junction; the continuous model by \cite{holden1995unidirection} is mathematically tractable, but the decoupling method based on an optimization problem is not directly associated with any physical merging or diverging rules.
In contrast, the third approach integrates the physical merging and diverging rules in CTM into the analytical framework of \cite{holden1995unidirection} is very promising for further studying network traffic dynamics.

Since the kinematic wave model of network traffic flow is a system of hyperbolic conservation laws on a network structure, solutions to the Riemann problem at a junction, in which all links carry constant initial conditions, but discontinuities can occur at the junction, are of physical, analytical, and numerical importance: physically, they can define physical merging, diverging, and other behavioral rules; analytically, a system of hyperbolic conservation laws is well-defined if and only if the Riemann problem is uniquely solved \cite{bressan2000convergence}; and numerically, they can be incorporated into the Godunov finite difference equations \cite{godunov1959}.
In this study, we present a Riemann solver of the network kinematic wave model using the solution framework of  \cite{jin2012_network}. 
Note that the Riemann solver is analytical in the sense of \cite{garavello2011conservation}, different from numerical ones as discussed in \cite{roe1981rieman}.  
The new solver is based on that in \cite{holden1995unidirection}: in the Riemann solutions, a stationary state arises on a link along with a shock or rarefaction wave, which is determined by the Riemann problem of the LWR model on the link with the initial and stationary states; the stationary state should be inside a feasible domain, such that the shock or rarefaction wave propagates backward on an upstream link and forward on a downstream link; and the constant in- or out-flux of a link equals the stationary flow-rate.    
The remaining piece in the Riemann solver is to introduce an entropy condition in the sense of \cite{holden1995unidirection} such that ``this condition gives a unique solution at least for Riemann initial data''.
But different from that in \cite{holden1995unidirection}, the Riemann solver in this study uses a discrete flux function, which is defined in terms of upstream demands and downstream supplies, as an entropy condition.
Such a discrete flux function is originally developed within the framework of CTM, models conflicts among merging and diverging traffic streams at the aggregate level, and is therefore a natural choice as the entropy condition to pick out unique physical solutions. 
To incorporate the new entropy condition, we enlarge the function space for weak solutions to the Riemann problem by introducing on each link an interior state, which is local and takes no space (of measure zero) right next to the junction. That is, the Riemann solution may be neither left- or right-continuous at jump discontinuities, and such a function space is still the same as that for traditional weak solutions. 
Then the entropy condition is introduced so that boundary fluxes through the junction be determined locally from  interior states with numerical CTM flux functions.
In addition, based on the equivalence between traffic density and the demand-supply pair, we solve the Riemann problem in the demand-supply space and show that stationary states exist and are unique. 
In this study, the local flux function is the global flux function derived in \cite{jin2012_network}.

The rest of the paper is organized as follows. In Section 2, we present the network kinematic wave model and review Holden and Risebro's Riemann solver and discrete CTM flux functions. In Section 3, we present a new Riemann solver. In Section 4, we solve the Riemann problem. In Section 5, we discuss some further properties of the Riemann solver. In Section 6, we present some concluding remarks.

\section{A system of hyperbolic conservation laws at a network junction and the Riemann problem}
In this study, we consider traffic dynamics in the junction network shown in \reff{general_junction}. We denote the set of upstream links by $A=\{1,\cdots,m\}$ and the set of downstream links by $B=\{m+1,\cdots,m+n\}$. 
For link $a\in A$, we introduce link coordinates $(a,x_a)$, where $x_a\in(-\infty, 0)$; for link $b\in B$, we introduce link coordinates $(b,x_b)$, where $x_b\in(0,\infty)$.
At a point $(a,x_a)$ ($a\in A\cup B$) and time $t$, we denote the total density, speed, and flow-rate by $k_a(x_a,t)$, $v_a(x_a,t)$, and $q_a(x_a,t)$, respectively. Hereafter we omit $(x_a,t)$ from these variables unless necessary.
Here we assume that vehicles of different paths, classes, or other attributes have the same characteristics, and that each link is homogeneous with a location-independent number of lanes, free-flow speeds, curvatures, slopes, and so on. Then we have the following fundamental diagram of flow-density and speed-density relations \cite{greenshields1935capacity}:
\bqs
v_a&=&V_a(k_a),\\
q_a&=&Q_a(k_a)\equiv k_a V_a(k_a).
\eqs 
Generally, $Q_a(k_a)$ is a unimodal function in $k_a$ and reaches its capacity, $C_a$, when traffic density equals the critical  density $k_{a,c}$. If traffic density $k_a$ is strictly smaller than, equal to, or strictly greater than the critical density $k_{a,c}$, then we call the traffic state as strictly under-critical (SUC), critical (C), or strictly over-critical (SOC), respectively. An under-critical state (UC) can be SUC or C, and an over-critical state (OC) can be SOC or C. Note that different links can have different fundamental diagrams.

From traffic conservation, we have the following system of LWR models in the network \cite{lighthill1955lwr,richards1956lwr}
\bqn \label{link-kw}
\pd{k_a}t+\pd{Q_a(k_a)}{x_a}&=&0, \quad a\in A\cup B.
\eqn
On a road link, it is well known that \refe{link-kw} admits weak solutions, among which the unique, physical solutions should satisfy so-called entropy conditions. For the LWR model, the traditional Lax's entropy condition \cite{lax1972shock} is consistent with vehicles' acceleration and deceleration behaviors \cite{ansorge1990entropy}.

On the network shown in \reff{general_junction}, \refe{link-kw} is a system of $m+n$ hyperbolic conservation laws, which can be understood as a semigroup based on the Godunov method and a Riemann solver \cite{bressan2000convergence}. Thus the main challenge of network kinematic wave theories is to develop a Riemann solver for the network, where hyperbolic conservation equations are coupled with each other at the junction due to merging and diverging conflicts among traffic streams. 
In the Riemann problem, we are interested in finding $k_a(x_a,t)$ at any $(x_a,t)$ from \refe{link-kw} with the following initial conditions with jump discontinuities at the junction:
\bsq\label{network-RP}
\bqn
 k_a(x_a,0)&=&k_a, \quad x_a<0, \qquad a\in A\\
k_b(x_b,0)&=&k_b, \quad x_b>0, \qquad b\in B
\eqn
\esq

For upstream link $a\in A$, we introduce a test function $\phi_a(x_a,t): (-\infty,0]\times[0,\infty)$, which is smooth with compact support on $(-\infty,0]\times[0,\infty)$; for downstream link $b\in B$, we introduce a test function $\phi_b(x_b,t): (0,\infty]\times[0,\infty)$, which is smooth with compact support on $(0,\infty]\times[0,\infty)$. These test functions are also smooth across the junction; i.e., $\phi_a(0,t)=\phi_b(0,t)$, and $\pd{\phi_a}{x_a}(0,t)=\pd{\phi_b}{x_b}(0,t)$ for $a\in A$ and $b\in B$. 
 Then we attempt to find the weak solutions to the Riemann problem of \refe{link-kw} with \refe{network-RP} in the following sense \cite{coclite2005network}:
\bqn
 \sum_{a\in A}\int_0^\infty \int_{-\infty}^0 (k_a \pd{\phi_a}{t}+Q_a(k_a)\pd{\phi_a}{x_a})d x_a dt+&&\nonumber\\\sum_{b\in B}\int_0^\infty \int_0^\infty (k_b \pd{\phi_b}{t}+Q_b(k_b)\pd{\phi_b}{x_b})d x_b dt&=&0. \label{weaksolution}
\eqn

\subsection{Holden and Risebro's Riemann solver}
 In \cite{holden1995unidirection}, the Riemann problem of \refe{link-kw} with \refe{network-RP} was decoupled into $m+n$ Riemann problems of the LWR model based on the following observations. (1) Due to similarity in Riemann solutions; i.e., $k_a(x_a,t)=\r_a(\frac{x_a}t)$, a stationary state, $k_a^*$, initiates at the junction and spreads on link $a$ ($a\in A\cup B$). 
 That is, for any $x_a\in(-\infty,0]$ ($a\in A$) or $x_a\in[0,\infty)$ ($a\in B$), $\lim_{t\to \infty}k_a(x_a,t)=k_a^*$, and the stationary state pervades the whole link after a long time. Note that it is possible that the stationary state is the same as the initial state.
 It is stationary in the sense that the boundary flux of link $a$ is constant and equals 
\bqs
 f_a=Q_a(k_a^*).
\eqs
  In addition, due to traffic conservation at the junction, we have 
\bqs
  \sum_{a\in A} f_a=\sum_{b\in B} f_b,
\eqs
which is also the Rankine-Hugoniot condition for the junction.
   (2) A shock or rarefaction wave can develop on each link and solves the Riemann problem for the corresponding LWR model with initial and stationary states as initial data pairs. That is, for upstream link $a$ ($a\in A$), the shock or rarefaction wave solves $\pd{k_a}t+\pd{Q_a(k_a)}{x_a}=0$ with $k_a(x_a,0)=\cas{{ll}k_a, & x_a<0\\k_a^*, & x_a>0}$;  for downstream link $b$ ($b\in B$), the shock or rarefaction wave solves $\pd{k_b}t+\pd{Q_b(k_b)}{x_b}=0$ with $k_b(x_b,0)=\cas{{ll}k_b^*, & x_b<0\\k_b, & x_b>0}$. 
Thus the Riemann problem is uniquely solved if and only if the stationary states are, and Holden and Risebro's Riemann solver is equivalent to finding the following mapping \cite{garavello2011conservation} 
\bqn
(k_1^*,\cdots,k_{m+n}^*)=\RS(k_1,\cdots,k_{m+n}). \label{HR-RS}
\eqn

Since stationary states can only propagate backward on upstream links and forward on downstream links, feasible regions of stationary states can be obtained by analyzing the $m+n$ Riemann problems of the LWR model.
Furthermore, as shown in \cite{holden1995unidirection}, there can exist multiple feasible solutions of stationary states. Thus an additional entropy condition has to be introduced to give a unique solution to the Riemann problem of \refe{link-kw} with \refe{network-RP}. In \cite{holden1995unidirection}, an entropy of a junction is defined by $E=\sum_{a=1}^{m+n} g(\frac{f_a}{C_a})$, and it was proved that there exists a well-defined Riemann solver when $g(\cdot)$ is a strictly concave function.
However, the entropy condition is not explicitly related to physical merging and diverging rules, and vehicles' pre-defined route choices are not considered.

In \cite{coclite2005network}, vehicles' pre-defined route choices are introduced by a matrix of turning proportions $\xi_{a\to b}$ ($\forall a\in A, b\in B$), where $\xi_{a\to b}\in[0,1]$, and 
\bqn
\sum_{b\in B} \xi_{a\to b}=1, \quad a\in A. \label{turningproportions}
\eqn
 Therefore, traffic conservation at the junction leads to $n$ equations: $\sum_{a\in A} f_a \xi_{a\to b} =f_b$ ($b\in B$). It was shown that the Riemann solver is well-defined for certain turning proportions. But this Riemann solver is not well-defined when $m>n$. 

In \cite{garavello2005source,garavello2006traffic,herty2006optimization,haut2007second,herty2008multicommodity}, more Riemann solvers along this line have been proposed, but the entropy conditions are not directly related to merging and diverging behaviors in these studies.

\subsection{Discrete CTM flux functions}

In \cite{daganzo1995ctm} and \cite{lebacque1996godunov}, the Godunov discrete form of the LWR model was extended to compute traffic flows through merging, diverging, and general junctions. In CTM, so-called traffic demand and supply, $d_a(x_a,t)$ and $s_a(x_a,t)$, are defined as functions of traffic density \cite{engquist1980difference,daganzo1995ctm,lebacque1996godunov}
\bsq \label{def:ds}
\bqn
d_a&=&D_a(k_a)\equiv Q_a(\min\{k_a,k_{a,c}\}),\\
s_a&=&S_a(k_a)\equiv Q_a(\max\{k_a,k_{a,c}\}),
\eqn
\esq
where traffic demand increases in total density $k_a$, and traffic supply decreases in total density $k_a$. Furthermore,  $q_a=\min\{d_a,s_a\}$, and $C_a=\max\{d_a,s_a\}$. In addition, since $d_a/s_a$ is a strictly increasing function of $k_a$, $k_a$ can be uniquely determined by $d_a/s_a$. We denote this function by
\bqn
k_a&=&K_a(d_a/s_a). \label{ds-k}
\eqn 
That is, the pair of demand and supply can uniquely determine the traffic state at a location and time. 

Within the framework of CTM, boundary fluxes through a junction at time $t$ are be calculated as follows: 
\bqn
(f_1^t,\cdots, f_{m+n}^t)=\FF^0(d_1^t,\cdots,d_m^t, s_{m+1}^t,\cdots, s_{m+n}^t), \label{fluxfunctions}
\eqn
where $d_1^t,\cdots, d_m^t$ are traffic demands in upstream cells next to the junction, $s_{m+1}^t,\cdots,s_{m+n}^t$ are traffic supplies in downstream cells next to the junction, and $f_1^t,\cdots,f_{m+n}^t$ are the out- and in-fluxes of all links.
The discrete flux function \refe{fluxfunctions} has been used to model traffic dynamics at various bottlenecks within the framework of CTM. For example, when vehicles have pre-defined route choices with turning proportions given in \refe{turningproportions}, a flux function was derived from the First-In-First-Out (FIFO) diverging and fair merging rules \cite{jin2004network}. Other CTM merging and diverging models can be found in \cite{jin2010merge,jin2010_diverge}.

Numerically, the flux functions \refe{fluxfunctions} can be incorporated into discrete traffic conservation equations and simulate traffic dynamics in a road network with given initial and boundary conditions.
Thus \refe{fluxfunctions} can be considered as approximate, numerical Riemann solvers \cite{leveque2002fvm}.
We can see that, using the concepts of demand and supply, it is rather straightforward to construct such flux functions for solving \refe{link-kw} numerically. While in many other systems of hyperbolic conservation laws, one has to solve the Riemann problem first and then obtain boundary fluxes as in the Godunov method. 
Physically, \refe{fluxfunctions} models how the right of way at the junction is allocated among competing traffic streams and represents drivers' macroscopic merging and diverging behaviors.
Some of these macroscopic diverging and merging rules have been verified through observations \cite{munoz2002diverge,lebacque2003intersection,ni2005merge,cassidy2005driver,bargera2010empirical}. 
 
\section{A new Riemann solver}
Since  discrete flux functions \refe{fluxfunctions} are both numerically and physically well-defined, it is reasonable to use it as an entropy condition to pick out unique, physical solution to the Riemann problem of \refe{link-kw} with \refe{network-RP}.
In \cite{jin2009sd,jin2010merge,jin2010_diverge,jin2012_network}, a new analytical framework was proposed by replacing the entropy condition in \cite{holden1995unidirection} by various CTM flux functions in \refe{fluxfunctions}.
In this section, we first describe the new framework and then focus on a new flux function derived in \cite{jin2012_network}.

\subsection{Riemann solvers with CTM flux functions as entropy conditions}\label{unifiedframework}
In order to use CTM flux functions \refe{fluxfunctions} as entropy conditions, we first enlarge the weak solution space by introducing a new interior state on each link in the Riemann solutions. That is, in the new weak solutions,  we assume that $k_a(0,t)=k_a^0$, but $\lim_{t\to\infty}k_a(x_a,t)=k_a^*$ for $x_a\in(-\infty,0)$ ($a\in A$) or $x_a\in(0,\infty)$ ($b\in B$).
Such an interior state is right adjacent to the junction and takes a space of measure zero.
Thus the new weak solutions to the Riemann problem still satisfy the traditional definition in \refe{weaksolution}.
In addition, they do not impact shock or rarefaction waves on all links, which are determined by stationary and initial states. 
 
Such interior states were first observed in numerical solutions of the Burgers equation \cite{vanleer1984upwind,bultelle1998shock} and the LWR model \cite{jin2003inhLWR} when stationary shock waves occur. Theoretically, if the LWR model $\pd{k}t+\pd{Q(k)}x=0$ on a line $x\in(-\infty,\infty)$ is solved by a zero-speed shock wave, then an interior state can be introduced at the stationary discontinuity without violating the definition of weak solutions or the Lax entropy condition. That is, $k(x,t)=\cas{{ll}k_L, &x<0\\k_M, &x=0\\k_R, &x>0}$ with $k_L\leq k_M \leq k_R$ and $Q(k_L)=Q(k_R)$ is a feasible solution with a flimsy interior state $k_M$ at $x=0$. 
In reality, such an interior state can be observed at the interface of a stationary shock wave, if a detector covers a part of the upstream traffic stream and a part of the downstream traffic stream. Depending on the relative location of the detector, such an interior state may not be unique. Thus the interior states can physically exist, and the enlarged function space of weak solutions is still well defined. In addition, in \cite{jin2010merge,jin2010_diverge,jin2012_network}, it was shown that the inclusion of interior states is necessary for Riemann solvers to be well-defined with some CTM flux functions.
Thus the introduction of interior states into weak solutions makes both physical and mathematical senses. From the viewpoint of traffic flow modeling, the new framework is much more powerful and flexible, since it allows many flux functions derived from driving rules as entropy conditions.

Since discrete flux functions \refe{fluxfunctions} are defined in demands and supplies, it is reasonable to use  demand and supply, instead of density, as state variables in the new Riemann solver; i.e., traffic condition at $(x_a,t)$ is determined by $U_a=(d_a,s_a)$. Thus a traffic state is UC if and only if $d_a\leq s_a=C_a$, or equivalently $U_a=(q_a,C_a)$; a traffic state is UC if and only if $s_a\leq d_a=C_a$, or equivalently $U_a=(C_a,q_a)$.

In the demand-supply space, we denote the initial and stationary traffic states on link $a$ by $(d_a,s_a)$ and $(d_a^*,s_a^*)$, respectively. Then we have the following observations in the Riemann solutions:
\bsq\label{feasiblestationary}
\bqn
f_a&=&\min\{d_a^*,s_a^*\}, \quad a\in A\cup B 
\eqn
In addition, we have the following lemma regarding the feasible regions of stationary states.
\begin{lemma} \cite{jin2012_network} In demand-supply space, the feasible regions of stationary states are given by
\bqn
U_a^*&\in& \B(U_a,\cdot)\equiv(d_a,C_a)\cup\{(C_a,s_a^*)|s_a^*<d_a\},\quad a\in A\\
U_b^*&\in& \F(\cdot, U_b)\equiv (C_b,s_b)\cup\{(d_b^*,C_b)|d_b^*<s_b\},\quad b\in B\\
\eqn
which lead to
\bqn
f_a&\leq &d_a,\quad  a\in A\\
f_b&\leq &s_b,\quad b\in B 
\eqn
\end{lemma}
\esq

Further, we denote the interior state on link $a$ by $U_a^0$ ($a\in A\cup B$).
Since interior states do not propagate into road links, the Riemann problems with stationary and interior states as initial data are solved by waves, whose speeds are positive on the upstream links and negative on the downstream links. Therefore, in supply-demand space, the feasible regions of interior states are given by
\bsq \label{feasibleinterior}
\bqn 
U_a^0\in \F(U_a^*,\cdot)\equiv\{(C_a, s_a^*)|s_a^*<d_a^*=C_a\}\cup\{(d_a^0,s_a^0)|s_a^0\geq d_a^*, d_a^*\leq s_a^*=C_a \}, \\
U_b^0\in \B(\cdot, U_b^*)\equiv\{(d_b^*, C_b)|d_b^*<s_b^*=C_b\}\cup\{(d_b^0,s_b^0)|d_b^0\geq s_b^*, s_b^*\leq d_b^*=C_b \}.
\eqn
\esq
where $a\in A$ and $b\in B$. When $U_a^0=U_a^*$, interior states do not exist.

From the feasible regions of stationary and interior states in \refe{feasiblestationary} and \refe{feasibleinterior}, we have the following corollary.
\begin{corollary} \label{cor:ss-is}
 For upstream link $a\in A$, $U_a^*$ is SOC if and only if $f_a<d_a$, and $U_a^*$ is UC if and only if $f_a=d_a$. In addition, if $U_a^*$ is SOC, then  $U_a^*=U_a^0=(C_a,f_a)$; if $U_a^*$ is UC, then $U_a^*=(d_a,C_a)$, and $U_a^0=(d_a^0,s_a^0)$ with $s_a^0\geq d_a$.
For downstream link $b\in B$, $U_b^*$ is SUC if and only if $f_b<s_b$, and $U_b^*$ is OC if and only if $f_b=s_b$. In addition, if $U_b^*$ is SUC, then $U_b^*=U_b^0=(f_b, C_b)$; if $U_b^*$ is OC, then $U_b^*=(C_b,s_b)$, and $U_b^0=(d_b^0, s_b^0)$ with $d_b^0 \geq s_b$.
\end{corollary}

Here we define four ratios for upstream link $a\in A$: the initial demand level $\delta_a=\frac{d_a}{C_a}$, the stationary demand level $\delta^*_a=\frac{d_a^*}{C_a}$, the interior demand level $\delta^0_a=\frac{d_a^0}{C_a}$, and the flux level $\tilde \delta_a=\frac{f_a}{C_a}$.
Then we have the following results.
\begin{corollary}\label{cor:demandlevels}
$U_a^*$ for $a\in A$ is SOC if and only if $\tilde \delta_a<\delta_a$; and $U_a^*$ is UC if and only if $\tilde \delta_a=\delta_a$. In addition, if $U_a^*$ is SOC, then 
\bqs
1=\delta^0_a=\delta^*_a\geq \delta_a>\tilde \delta_a;
\eqs
if $U_a^*$ is UC, then
\bqs
\delta^*_a=\delta_a=\tilde \delta_a.
\eqs
\end{corollary}
Note that the relationship between $\delta^0_a$ and other demand levels when $U_a^*$ is UC is to be determined.

Similarly we define four ratios for downstream link $b\in B$: the initial supply level $\sigma_b=\frac{s_b}{C_b}$, the stationary supply level $\sigma^*_b=\frac{s_b^*}{C_b}$, the interior supply level $\sigma^0_b=\frac{s_b^0}{C_b}$, and the flux level $\tilde \sigma_b=\frac{f_b}{C_b}$.
Then we have the following results.
\begin{corollary}\label{cor:supplylevels}
$U_b^*$ for $b\in B$ is SUC if and only if $\tilde \sigma_b<\sigma_b$; and $U_b^*$ is OC if and only if $\tilde \sigma_b=\sigma_b$. In addition, if $U_b^*$ is SUC, then
\bqs
1=\sigma^0_b=\sigma^*_b\geq \sigma_b>\tilde \sigma_b;
\eqs
if $U_b^*$ is OC, then
\bqs
\sigma^*_b=\sigma_b=\tilde \sigma_b.
\eqs
\end{corollary}
Note that the relationship between $\sigma^0_b$ and other supply levels when $U_b^*$ is OC is to be determined.

In the new Riemann solvers, various CTM flux functions in \refe{fluxfunctions} are used to determine fluxes locally from interior states as follows:
\bqn
(f_1,\cdots,f_{m+n})=\FF^0(d_1^0,\cdots,d_m^0,s_{m+1}^0,\cdots,s_{m+n}^0). \label{interiorflux}
\eqn 
From the corollaries above, we can see that, given initial upstream demands, $d_a$ ($a\in A$), and downstream supplies, $s_b$ ($b\in B$), if one can find all boundary fluxes, $f_a$ and $f_b$, then the stationary states can be uniquely determined. But the interior states may not be uniquely determined. Thus the Riemann problem is uniquely solved in the sense that all stationary states and waves on all links are uniquely solved.
Therefore, the Riemann problem is uniquely solved if and only if, from \refe{feasiblestationary}, \refe{feasibleinterior}, and \refe{interiorflux}, we can find a unique flux function that maps initial conditions into boundary fluxes
\bqn
(f_1,\cdots,f_{m+n})&=&\FF(d_1,\cdots,d_m,s_{m+1},\cdots,s_{m+n}). \label{initialflux}
\eqn
Hereafter, we refer to $\FF^0(\cdots)$ in \refe{interiorflux} as local or discrete flux functions and $\FF(\cdots)$ in \refe{initialflux} as global or continuous flux functions. Obviously, $\FF(\cdots)$ is the Godunov flux function, since it is derived from Riemann solutions.
 
Then from Corollary \ref{cor:ss-is}, we obtain the corresponding demand-supply Riemann solver, which maps upstream demands and downstream supplies into stationary states:
\bqn
(U_1^*,\cdots,U_{m+n}^*)&=&\RS'(d_1,\cdots,d_m,s_{m+1},\cdots,s_{m+n}). \label{ds-RS}
\eqn
Since both demand and supply are many-to-one functions of density, it suggests that different initial densities could lead to the same stationary states, if and only if the corresponding upstream demands and downstream supplies are the same. From $\RS'(\cdots)$ we can obtain the traditional Riemann solver $\RS(\cdots)$ in \refe{HR-RS} by first converting $(k_1,\cdots,k_{m+n})$ into $(d_1,\cdots,d_m,s_{m+1},\cdots, s_{m+n})$, then applying \refe{ds-RS} to find stationary demands and supplies, and finally converting $(U_1^*,\cdots,U_{m+n}^*)$ into $(k_1^*,\cdots,k_{m+n}^*)$ using \refe{ds-k}. However, from the traditional Riemann solver one may not be able to obtain a demand-supply Riemann solver.
It can be seen that a Riemann solver in \refe{ds-RS} satisfies the consistency condition \cite{garavello2011conservation}, if and only if $\FF(d_1,\cdots,d_m,s_{m+1},\cdots,s_{m+n})=\FF(d_1^*,\cdots,d_m^*,s_{m+1}^*,\cdots,s_{m+n}^*)$, which is equivalent to that $(U_1^*,\cdots,U_{m+n}^*)=\RS'(d_1^*,\cdots,d_m^*,s_{m+1}^*,\cdots,s_{m+n}^*)$. 
In addition, we call the local or discrete flux function $\FF^0(\cdots)$ invariant if $\FF^0(\cdots)=\FF(\cdots)$.
 
\subsection{A new local flux function}
In this study, we use the Godunov flux function in \cite{jin2012_network} as a new local flux function, \refe{interiorflux}. Note that the discrete flux function in \cite{jin2012_network} is consistent with fair merging and first-in-first-out diverging rules but different from the Godunov flux function.
Here all vehicles have predefined route choices, and the turning proportions, $\xi_{a\to b}$, are given in \refe{turningproportions}. The discrete flux function $\FF^0(\cdots)$ is defined as follows:
\ben
\item The out-flux of upstream link $a\in  A $ is given by
\bsq\label{entropy2}
\bqn
f_a&=&\min\{d_a^0,\theta^0 C_a\}, \label{entropy:upstream}
\eqn
where the interior critical demand level, $\theta^0$, is defined as follows:
\bqn
\theta^0&=&\min\left\{\max_{a\in A} \frac{d_a^0}{C_a},\min_{b\in  B   } \max_{A_1\subseteq A} \frac{s^0_b-\sum_{\a\in A\setminus A_1} d_\a^0 \xi_{\a\to b}}{\sum_{a\in A_1} C_a \xi_{a\to b}}\right\} , \label{deftheta}
\eqn
where $A_1$ is not empty. 
\item The in-flux of downstream link $b\in  B   $ is given by
\bqn
f_b&=&\sum_{a\in  A  } f_a \xi_{a\to b}. \label{entropy:downstream}
\eqn
\esq
\een
In terms of demand and supply levels, $\FF^0(\cdots)$ in \refe{entropy2} can be re-written as
\bsq \label{newfluxfunction}
\bqn 
\tilde \delta_a&=&\min\{\delta_a^0,\theta^0\},\\
\theta^0&=&\min\left\{\max_{a\in A} \delta_a^0,\min_{b\in  B   } \max_{A_1\subseteq A} \frac{C_b \sigma^0_b-\sum_{\a\in A\setminus A_1} C_{\a\to b} \delta_\a^0}{\sum_{a\in A_1} C_{a\to b}}\right\} ,\label{newtheta}\\
C_b \tilde \sigma_b&=& \sum_{a\in A} C_{a\to b} \tilde \delta_a,
\eqn
\esq
where $C_{a\to b}=C_a\xi_{a\to b}$ for $a\in A$ and $b\in B$.

\subsection{Average demand levels}
In this subsection, we discuss  properties of $\theta^0$ in \refe{newtheta}.
For $a\in A$ and $b\in B$, we define the demand and supply levels by $\mu_a\in [0,1]$ and $\nu_b\in [0,1]$, respectively. Further, we denote $\pi_b=C_b \nu_b - \sum_{a\in A}C_{a\to b} \mu_a$. \footnote{In this study we do not consider situations when $\xi_{a\to b}=0$.} We denote the vector of $\mu_a$ for $a\in A$ by $\bfmu$, and the vector of $\nu_b$ for $b\in B$ by $\bfnu$.
We define the average demand level of set $A_1$ for link $b$ by
\bqn
\gamma_b(A_1)&=&\frac{C_b \nu_b-\sum_{\a\in A\setminus A_1} C_{\a\to b} \mu_\a}{\sum_{a\in A_1} C_{a\to b}}=\frac{\pi_b+\sum_{a\in A_1} C_{a\to b} \mu_a}{\sum_{a\in A_1} C_{a\to b}}, \label{averagedl}
\eqn
where $A_1\subseteq A$ is non-empty.
\begin{lemma}\label{lemma:gamma}
For $\a\in A_1$, and $A_2\equiv A_1 \setminus \{\a\}\neq \emptyset$, if $\mu_\a>,=,<\gamma_b(A_1)$, then $\gamma_b(A_2)<,=,>\gamma_b(A_1)$, and $\mu_\a>,=,<\gamma_b(A_2)$, respectively.
For $\a \notin A_1$, and $A_2\equiv A_1 \cup \{\a\}$, if $\mu_\a<,=,>\gamma_b(A_1)$,  then $\gamma_b(A_2)<,=,>\gamma_b(A_1)$, and $\mu_\a<,=,>\gamma_b(A_2)$, respectively.
That is, if we remove a link with a larger demand level, then the average demand level decreases; if we add a link with a larger demand level, then the average demand level increases.
\end{lemma}

{\em Proof}.  For $\a\in A_1$, and $A_2\equiv A_1 \setminus \{\a\}\neq \emptyset$, we have
\bqs
\gamma_b(A_1)&=&\frac{\sum_{a\in A_2} C_{a\to b}\mu_a+C_{\a\to b}\mu_\a+\pi_b}{\sum_{a\in A_2} C_{a\to b}+C_{\a\to b}},
\eqs
which leads to
\bqs
\gamma_b(A_2)&=&\frac{\sum_{a\in A_2} C_{a\to b}\mu_a+\pi_b}{\sum_{a\in A_2} C_{a\to b}}=\gamma_b(A_1)+\frac{(\gamma_b(A_1) -\mu_\a) C_{\a\to b}}{\sum_{a\in A_2}C_{a\to b}},
\eqs
Thus, if $\mu_\a>,=,<\gamma_b(A_1)$, then  $\gamma_b(A_2)<,=,>\gamma_b(A_1)$, and $\mu_\a>,=,<\gamma_b(A_2)$, respectively.

For $\a \notin A_1$, and $A_2\equiv A_1 \cup \{\a\}$, we have
\bqs
\gamma_b(A_1)&=&\frac{\sum_{a\in A_2} C_{a\to b}\mu_a-C_{\a\to b}\mu_\a+\pi_b}{\sum_{a\in A_2} C_{a\to b}-C_{\a\to b}},
\eqs
which leads to
\bqs
\gamma_b(A_2)&=&\frac{\sum_{a\in A_2} C_{a\to b} \mu_a+\pi_b}{\sum_{a\in A_2} C_{a\to b}}= \gamma_b(A_1)-\frac{(\gamma_b(A_1)-\mu_\a) C_{\a\to b}}{\sum_{a\in A_2}C_{a\to b}}.
\eqs
Thus, if $\mu_\a<,=,>\gamma_b(A_1)$, $\gamma_b(A_2)<,=,>\gamma_b(A_1)$, respectively.
Furthermore, 
\bqs
(\gamma_b(A_2)-\mu_\a) \sum_{a\in A_2} C_{a\to b}&=&\sum_{a\in A_1} C_{a\to b}\mu_a +C_{\a\to b}\mu_\a+\pi_b- \mu_\a \sum_{a\in A_1} C_{a\to b}-C_{\a\to b}\mu_\a
\\&=&\sum_{a\in A_1} C_{a\to b}\mu_a +\pi_b- \mu_\a \sum_{a\in A_1} C_{a\to b}=(\gamma_b(A_1)-\mu_\a) \sum_{a\in A_1} C_{a\to b}.
\eqs
Thus, if $\mu_\a<,=,>\gamma_b(A_1)$, $\mu_\a<,=,>\gamma_b(A_2)$, respectively.
 \eop

We denote the maximum average demand level for link $b$ by $\Gamma_b=\max_{A_1\subseteq A, A_1\neq \emptyset} \gamma_b(A_1)$, which has the following properties.
\begin{lemma} \label{lemma:averagedl} $\gamma_b(A_1)$ and $\Gamma_b$ have the following properties.   
\ben
\item When $\pi_b>0$, then $\gamma_b(\{a\})>\mu_a$ for $a\in A$, and $\Gamma_b=\max_{a\in A} \gamma_b(\{a\})>\max_{a\in A} \mu_a$.
\item When $\pi_b=0$, then $\gamma_b(\{a\})=\mu_a$ for $a\in A$, and $\Gamma_b=\max_{a\in A} \gamma_b(\{a\})=\max_{a\in A} \mu_a$.
\item When $\pi_b<0$, then $\gamma_b(\{a\})<\mu_a$ for $a\in A$, and there exists a unique $A_1\neq \emptyset$, such that $\Gamma_b=\gamma_b(A_1)$, $\mu_a>\Gamma_b\geq \mu_\a$ for $a\in A_1$ and $\a \in A\setminus A_1$.
\een
\end{lemma}
{\em Proof}. 
\ben
\item When $\pi_b>0$, from \refe{averagedl}, we can have $\gamma_b(\{a\})>\mu_a$ for $a\in A$. Then from Lemma \ref{lemma:gamma} we have 
$\min_{a\in A_1} \mu_a < \gamma_b(A_1)\leq\max_{a\in A_1}\gamma_b(\{a\})$. Thus $\Gamma_b=\max_{a\in A} \gamma_b(\{a\})>\max_{a\in A} \mu_a$. Note that $\max_{a\in A} \gamma_b(\{a\})$ and $\max_{a\in A} \mu_a$ may attain their maxima for different $a$.
\item When $\pi_b=0$, from \refe{averagedl}, we can have $\gamma_b(\{a\})=\mu_a$ for $a\in A$. Then from Lemma \ref{lemma:gamma} we have 
$\min_{a\in A_1} \mu_a\leq \gamma_b(A_1)\leq\max_{a\in A_1} \mu_a$. Thus $\Gamma_b=\max_{a\in A} \gamma_b(\{a\})=\max_{a\in A} \mu_a$. 
\item When $\pi_b<0$, from \refe{averagedl}, we can have $\gamma_b(\{a\})<\mu_a$ for $a\in A$. Since $A$ has a finite number of subsets $A_1$,  $\Gamma_b$ exists and we can find $A_1$, such that $\Gamma_b=\gamma_b(A_1)$, and the value of $\Gamma_b$ is unique. First, for $\a\notin A_1$, $\mu_\a \leq \gamma_b(A_1)$, since, otherwise, from Lemma \ref{lemma:gamma} we have $\gamma_b(A_1)<\gamma_b(A_1\cup \{\a\})$, which contradicts that $\Gamma_b=\gamma_b(A_1)$. 
Second, there exists at least one $a\in A_1$, such that $\mu_a>\gamma_b(A_1)$, since, otherwise, from Lemma \ref{lemma:gamma} we have $\gamma_b(\{a\})=\mu_a$ for all $a\in A_1$, which contradicts that $\gamma_b(\{a\})<\mu_a$.
Third, for all $a\in A_1$, $\mu_a\geq \gamma_b(A_1)$, since, otherwise, from Lemma \ref{lemma:gamma} we have $\gamma_b(A_1)<\gamma_b(A_1\setminus \{a\})$, which contradicts that $\Gamma_b=\gamma_b(A_1)$.  Assume that $a\in A_1^*\subset A_1$ and $\mu_a>\gamma_b(A_1)$, then $\mu_\a=\gamma_b(A_1)$ for $\a \in A_1\setminus A_1^*$, and $\gamma_b(A_1^*)=\gamma_b(A_1)$. Since $A_1^*\neq \emptyset$ is unique, there exists a unique $A_1=A_1^*\neq \emptyset$, such that $\Gamma_b=\gamma_b(A_1)$, $\mu_a>\Gamma_b\geq \mu_\a$ for $a\in A_1$ and $\a \in A\setminus A_1$.     
\een
  \eop

Assuming that $\mu_a$ ($a\in A$) are in a decreasing order; i.e., $\mu_1\geq \cdots \geq \mu_m$, we define the following average demand level of the first $l$ upstream links:
\bqs
\gamma_b(l)&=&\gamma_b(\{1,\cdots,l\}).
\eqs
Then we have the following corollary.
\begin{corollary} $\gamma_b(l)$ and $\Gamma_b$ have the following properties:
\ben
\item When $\pi_b>0$, $\gamma_b(l)>\mu_l$ for $l\in A$, $\gamma_b(1)>\cdots>\gamma_b(m)$, and $\Gamma_b=\max_{a\in A} \gamma_b(\{a\})\geq \gamma_b(1)>\mu_1$.
\item When $\pi_b=0$, $\gamma_b(l)\geq \mu_l$ for $l\in A$, $\gamma_b(1)\geq\cdots\geq\gamma_b(m)$, and $\Gamma_b=\mu_1$.
\item When $\pi_b<0$,  there exists a unique $l^*\in \{1,\cdots,m\}$, such that $\gamma_b(1)<\cdots <\gamma_b(l^*)\geq \gamma_b(l^*+1)\geq \cdots \geq \gamma_b(m)$, $\gamma_b(l)<\mu_l$ for $l=1,\cdots,l^*$, and $\gamma_b(l)\geq \mu_l$ for $l=l^*+1,\cdots,m$. That is, for every $b$, there exists a unique $l^*_b$ solving $\min l$, such that $\gamma_b(l)\geq \gamma_b(l+1)$. In addition, $l^*_b=1$ and $\gamma_b(1) \geq \mu_1$ if and only if $\pi_b\geq 0$. 
\een
\end{corollary}
{\em Proof}. When $\pi_b\geq 0$, the results are obvious. Here we will focus on $\pi_b<0$. We denote $l^*$ as the number of links in $A_1$, where $\Gamma_b=\gamma_b(A_1)$ and $\mu_a>\Gamma_b\geq \mu_\a$ for $a\in A_1$ and $\a\in A\setminus A_1$. Then $l^*$ is unique, and $\mu_{l^*}>\Gamma_b=\gamma_b(l^*)\geq \mu_b(l^*+1)$. From Lemma \ref{lemma:gamma}, we have $\gamma_b(1)<\cdots <\gamma_b(l^*)\geq \gamma_b(l^*+1)\geq \cdots \geq \gamma_b(m)$, $\gamma_b(l)<\mu_l$ for $l=1,\cdots,l^*$, and $\gamma_b(l)\geq \mu_l$ for $l=l^*+1,\cdots,m$.
 \eop

We define the critical demand level by
\bqn
g(\bfmu,\bfnu)&=&\min_b\max_{A_1}\frac{C_b \nu_b-\sum_{\a\in A\setminus A_1} C_{\a\to b} \mu_\a}{\sum_{a\in A_1} C_{a\to b}}\nonumber\\
&=&\min_b\max_{A_1}\frac{\pi_b+\sum_{a\in A_1} C_{a\to b} \mu_a}{\sum_{a\in A_1} C_{a\to b}} \label{generaltheta}
\eqn
Since $g(\bfmu,\bfnu)=\min_b \Gamma_b$, we have the following theorem regarding $g(\bfmu,\bfnu)$.
\begin{theorem} \label{thm:cdl}
For $C_{a\to b}>0$, $C_b>0$, $\mu_a\in [0,1]$, and $\nu_b\in [0,1]$ ($a\in A$ and $b\in B$), $g(\bfmu, \bfnu)$ defined in \refe{generaltheta} has the following properties:
\ben
\item $g(\bfmu, \bfnu)\geq \max_a \mu_a$ if and only if $\min_{b\in B}\pi_b\geq 0$.
\item If and only if $\min_{b\in B}\pi_b< 0$, there exists a unique non-empty $A_*\subseteq A$ such that
\bqn
g(\bfmu,\bfnu)&=&\min_b \max_{A_1} \gamma_b(A_1)=\min_b\gamma_b(A_*)=\max_{A_1} \min_b \gamma_b(A_1), \label{minmax}
\eqn
and
\bqs
\min_{a\in A_*} \mu_a > g(\bfmu,\bfnu)\geq \max_{\a \in A\setminus A_*} \mu_a.
\eqs
\een
\end{theorem}
{\em Proof}. When $\min_b\pi_b\geq 0$, from Lemma \ref{lemma:averagedl} we have  $g(\bfmu, \bfnu)=\min_b \Gamma_b \geq \max_a \mu_a$.

When $\min_b\pi_b<0$, we denote $A_b=\{a\in A| \mu_a> \Gamma_b\}$. From Lemma \ref{lemma:averagedl}, $A_b=\emptyset$ if and only if $\pi_b\geq 0$. We denote $A_*=\cup_{b\in B} A_b \neq \emptyset$, since $\min_b \pi_b<0$. In addition, for any $a\in A_*$, $\mu_a>\min_b \Gamma_b=g(\bfmu,\bfnu)$; and for any $\a\in A\setminus A_*$, $\mu_a\leq \min_b \Gamma_b=g(\bfmu,\bfnu)$. In addition, $g(\bfmu,\bfnu)=\min_b \gamma_b(A_*)$, and $\min_{a\in A_*} \mu_a >\min_b \gamma_b(A_*)\geq \max_{\a \in A\setminus A_*} \mu_a$.

Since $\min_b \pi_b<0$, we have $\min_b \gamma_b(\{a\})<\mu_a$ for $a\in A$. Since $A$ has a finite number of subsets $A_1$, thus we can find $A_1^*$ such that $\max_{A_1} \min_b\gamma_b(A_1)=\min_b\gamma_b(A_1^*)$, and the maximum value is unique. First, for any $\a\notin A_1^*$ $\mu_\a\leq \min_b\gamma_b(A_1^*)$, since, otherwise, from Lemma \ref{lemma:gamma} $\min_b\gamma_b(A_1^*\cup\{\a\})> \min_b\gamma_b(A_1^*)$.
Thus $a\in A_1^*$ if $\mu_a> \min_b\gamma_b(A_1^*)$.
Second, there exists at least one $a\in A_1^*$, such that $\mu_a> \min_b\gamma_b(A_1^*)$, since, otherwise, $\mu_a \leq \min_b\gamma_b(A_1^*)$ for all $a\in A_1^*$, and from Lemma \ref{lemma:gamma} $\mu_a \leq \min_b\gamma_b(\{a\})$, which is not possible. 
Third, if $a\in A_1^*$, then $\mu_a\geq \min_b\gamma_b(A_1^*)$, since, otherwise, from Lemma \ref{lemma:gamma} $\min_b\gamma_b(A_1^*\setminus\{a\})> \min_b\gamma_b(A_1^*)$. 
Therefore, without loss of generality, we can remove all $a$ from $A_1^*$ if $\mu_a=\min_b\gamma_b(A_1^*)$, and the value of $\min_b\gamma_b(A_1^*)$ does not change. Then $A_1^*=\{a\in A| \mu_a>\min_b\gamma_b(A_1^*)\}$ is unique; that is, there exists a unique $A_1^*$ such that $\min_{a\in A_1^*} \mu_a>\min_b\gamma_b(A_1^*) \geq \max_{\a\in A\setminus A_1^*} \mu_\a$.

Therefore, $A_*=A_1^*$, and \refe{minmax} is proved.  \eop

From Theorem \ref{thm:cdl} we have the following corollary.
\begin{corollary} If $\mu_a$ is decreasingly ordered, we have the following results on $g(\bfmu,\bfnu)$:
\ben
\item If and only if $\min_b\pi_b\geq 0$, $g(\bfmu,\bfnu)\geq \mu_1$.
\item If and only if $\min_b\pi_b< 0$, there exists a unique $l^*\in\{1,\cdots,m\}$, such that
\bqn
g(\bfmu,\bfnu)&=&\min_b\max_l \gamma_b(l) =\min_b \gamma_b(l^*)=\max_l \min_b \gamma_b(l),
\eqn
and
\bqn
\mu_{l^*}>g(\bfmu,\bfnu)\geq \mu_{l^*+1}.
\eqn
\een
\end{corollary}

\begin{figure}
\bc
\includegraphics[width=4in]{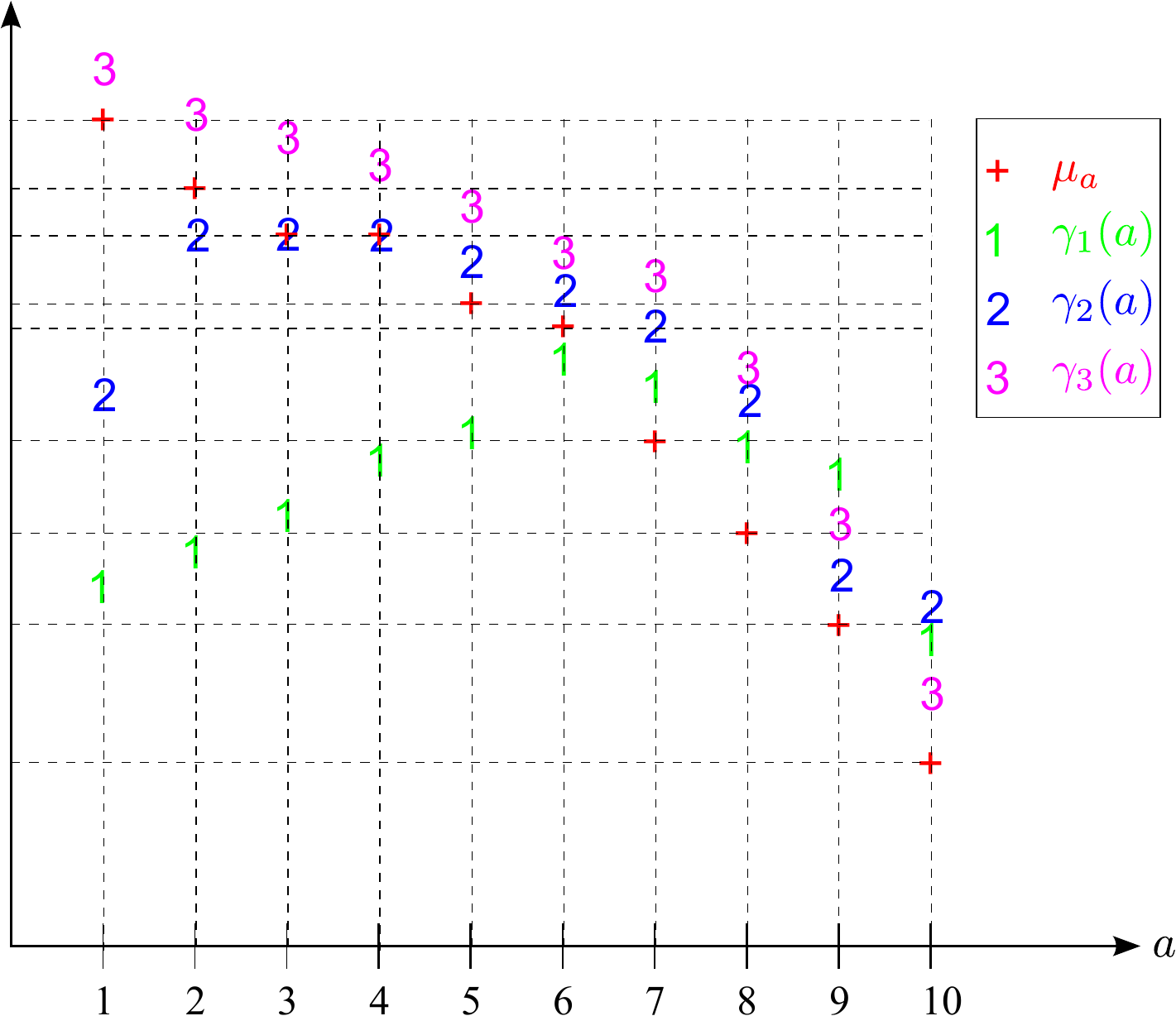}\caption{The pattern of $\gamma_b(a)$}\label{gamma_pattern}
\ec
\end{figure}

We show an example in \reff{gamma_pattern}, where $\mu_1\geq \cdots \geq \mu_{10}$. 
For $b=1$, $\gamma_1(1)<\cdots<\gamma_1(6)>\gamma_1(7)>\cdots>\gamma_1(10)$, $\gamma_2(a)<\mu_a$ for $a=1,\cdots,6$, and $\gamma_2(a)\geq \mu_a$ for $a=7,\cdots,10$. 
For $b=2$, $\gamma_2(1)<\gamma_2(2)=\gamma_2(3)=\gamma_2(4)>\cdots>\gamma_2(10)$, $\gamma_2(a)<\mu_a$ for $a=1,2$, and $\gamma_2(a)\geq \mu_a$ for $a=3,\cdots,10$.
For $b=3$ $\gamma_3(1)> \cdots > \gamma_3(10)$, and $\gamma_3(a)\geq \mu_a$ for $a=1,\cdots,10$.  
Then we can verify the lemmas above and find that $l^*=6$, such that $\mu_6>g(\bfmu,\bfnu)=\gamma_1(6)\geq \mu_7$.

Since $\theta^0=\min\{\max_{a\in A} \delta_a^0, g(\bfdelta^0,\bfsigma^0)\}$, where $\bfdelta^0=(\delta_a^0)_{a\in A}$ and $\bfsigma^0=(\sigma_b^0)_{b\in B}$, we have the following observations from Theorem \ref{thm:cdl}:
\ben
\item When $C_b \sigma_b^0\geq \sum_{a\in A}C_{a\to b}\delta_a^0$ for all $b\in B$, then $g(\bfdelta^0,\bfsigma^0)\geq \max_{a\in A} \delta_a^0$, and $\theta^0=\max_{a\in A} \delta_a^0$.
\item When $C_b \sigma_b^0< \sum_{a\in A}C_{a\to b}\delta_a^0$ for some $b\in B$, then there exists a unique, non-empty $A_*$ such that $\theta^0=g(\bfdelta^0,\bfsigma^0)=\min_b \gamma_b(A_*)$, and $\min_{a\in A_*} \delta_a^0>\theta^0\geq \max_{\a\in A\setminus A_*}\delta_\a^0$.
\een 
Thus $\theta^0$ is well-defined, bounded between $0$ and $\max_{a\in A} \delta_a^0$, and continuous in $(\delta_1^0,\cdots,\delta_m^0,\sigma_{m+1}^0,\cdots,\sigma_{m+n}^0)$.

\section{Solutions to the Riemann problem}

From the preceding section, we can see that the demand and supply levels satisfy the following conditions:
\bsq\label{finalrelations}
\bqn
&\{1=\delta^0_a=\delta^*_a\geq \delta_a>\tilde \delta_a\}\m{ or } \{\delta_a^*=\delta_a=\tilde \delta_a \m{ and } s_a^0\geq d_a\}& \\
&\{1=\sigma^0_b=\sigma^*_b\geq \sigma_b>\tilde \sigma_b\}\m{ or } \{\sigma_b^*=\sigma_b=\tilde \sigma_b \m{ and } s_b^0\geq d_b\}& \\
&\tilde \delta_a=\min\{\delta_a^0,\theta^0\}& \\
&\theta^0=\max_{a\in A}\delta_a^0\m{ or } (\delta_a^0)_{a\in A_*\neq \emptyset}>\theta^0=\min_{b\in B} \gamma_b^0(A_*) \geq (\delta_a^0)_{a\in A\setminus A_*}& \\
&C_b \tilde \sigma_b=\sum_{a\in A} C_{a\to b} \tilde \delta_a&
\eqn
\esq
where $\gamma_b^0(A_*)=\frac{C_b \sigma^0_b-\sum_{\a\in A\setminus A_*} C_{\a\to b} \delta_\a^0}{\sum_{a\in A_*} C_{a\to b}}$.
In this section, we attempt to solve $\tilde \delta_a$ ($a\in A$) in $\delta_a$ and $\sigma_b$; i.e., a mapping from upstream demand levels and downstream supply levels to upstream flux levels:
\bqn
(\tilde \delta_1,\cdots,\tilde \delta_m)&=&\FF'(\delta_1,\cdots,\delta_m,\sigma_{m+1},\cdots,\sigma_{m+n}).
\eqn
Then from the definition of $\tilde \delta_a$ and \refe{entropy:downstream} we can have the flux function $(f_1,\cdots,f_{m+n})=\FF(d_1,\cdots,d_m,s_{m+1},\cdots,s_{m+n})$ and the corresponding Riemann solver.

\subsection{Further properties of demand and supply levels}
\begin{figure}
\bc
\includegraphics[width=4in]{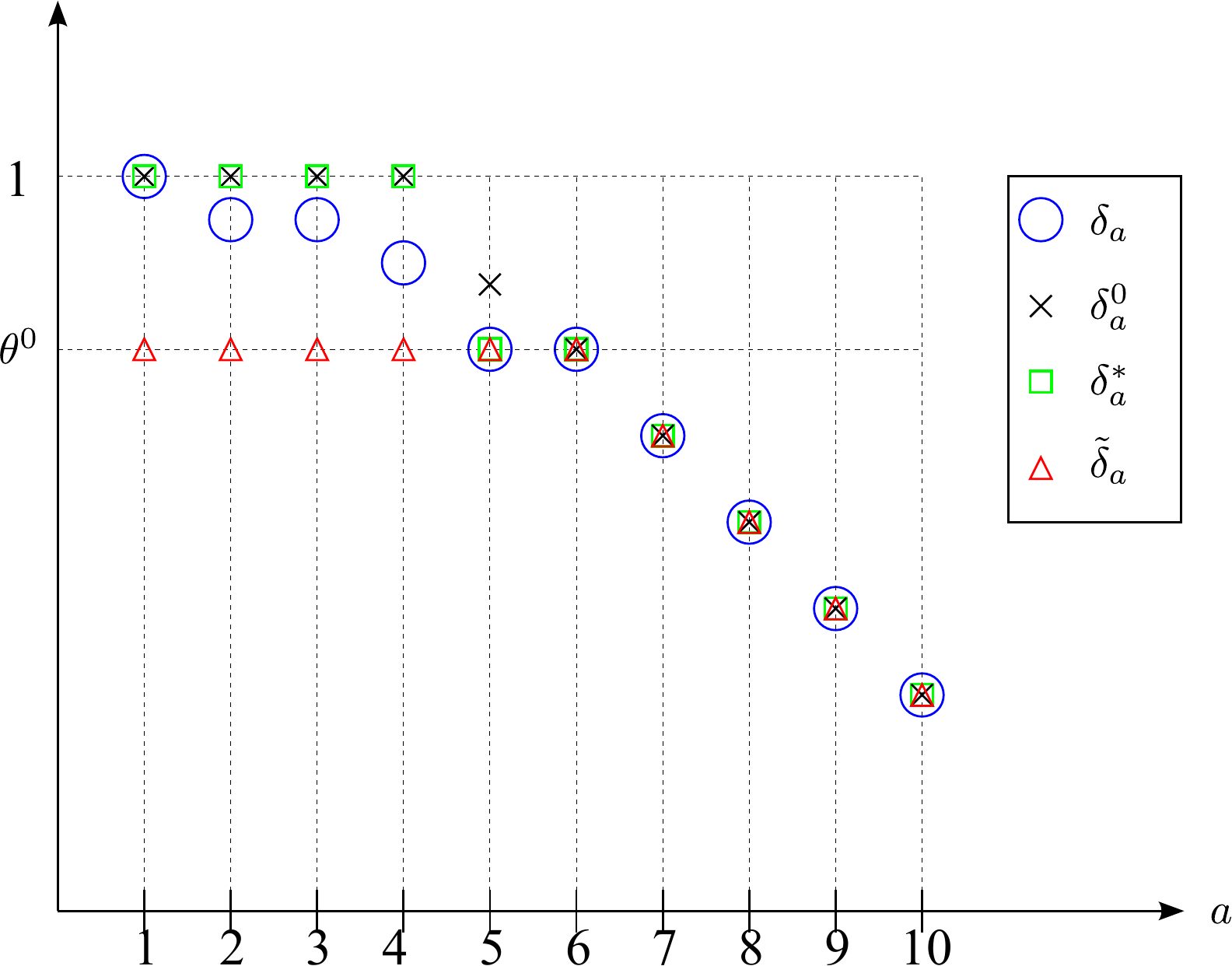}\caption{The relationships between $\delta_a$, $\delta^*_a$, $\delta^0_a$, $\tilde \delta_a$, and $\theta^0$}\label{demand_levels}
\ec
\end{figure}

From Corollary \ref{cor:demandlevels} and \refe{finalrelations} we have the following lemma.
\begin{lemma} \label{lemma:updemand}
$U_a^*$ is SOC if and only if $\theta^0=\tilde \delta_a<\delta_a$, and $U_a^*$ is UC if and only if $\theta^0\geq \delta_a$. In addition, if $\theta^0<\delta_a$, then $U_a^*=U_a^0=(C_a,f_a)$ is SOC, and
\bqs
1=\delta^0_a=\delta^*_a\geq \delta_a>\tilde \delta_a=\theta^0;
\eqs
if $\theta^0>\delta_a$, then $U_a^*=U_a^0=(d_a,C_a)$ is SUC, and 
\bqs
\delta^0_a= \delta^*_a=\delta_a=\tilde \delta_a<\theta^0;
\eqs
if $\theta^0=\delta_a$, then $U_a^*=(d_a,C_a)$ is UC, $U_a^0=(d_a^0,s_a^0)$ with $d_a^0\geq d_a$ and $s_a^0\geq d_a$, and 
\bqs
\delta^0_a\geq \delta^*_a=\delta_a=\tilde \delta_a=\theta^0.
\eqs
In the third case, $U_a^0=U_a^*$ if and only if $\delta^0_a=\theta^0$, and the interior state $U_a^0$ is different from the stationary state $U_a^*$ when $\delta^0_a>\theta^0$.
\end{lemma}
An example of the relationships between $\delta_a$, $\delta_a^*$, $\delta^0_a$, and $\theta^0$ is shown in \reff{demand_levels}, in which $U_a^*=U_a^0$ is SOC for $a=1,\cdots,4$, $U_a^*=U_a^0$ is UC for $a=6,\cdots,10$, and $U_5^*$ is UC but $U_5^*\neq U_5^0$.

We define $A_*^0=\{a\in A| \delta^0_a>\theta^0\}$. From Lemma \ref{lemma:updemand}, $U_a^*$ is either SOC or UC with $U_a^*\neq U_a^0$ for $a\in A_*^0$; and $d_\a^0=d_\a$ for $\a \in A\setminus A_*^0$. 
From Theorem \ref{thm:cdl}, we can see that $\min_{b\in B}\pi_b^0=\min_{b\in B}s_b^0-\sum_{a=1}^m C_{a\to b}\delta^0_a \geq 0$ iff $A_*^0$ is empty.
If $A_*^0$ is non-empty, then $(\delta_a^0)_{a\in A_*^0} >\theta^0=\min_{b\in B}\gamma^0_b(A_*^0)\geq (\delta_\a^0)_{\a\in A\setminus A_*^0}$.
Then we have the following lemma. 
\begin{lemma} \label{lemma:downsupply}
For downstream link $b$, if $A_*^0\neq \emptyset$ and $\gamma^0_b(A_*^0)=\theta^0$ or if $A_*^0=\emptyset$ and $\pi_b^0=0$, then $U_b^*=U_b^0=(C_b,s_b)$ is OC, and
\bqs
\sigma^0_b= \sigma^*_b=\sigma_b=\tilde \sigma_b;
\eqs
When $A_*^0\neq \emptyset$ and $\gamma^0_b(A_*^0)>\theta^0$ or when $A_*^0=\emptyset$ and $\pi_b^0>0$, there can be two types of solutions: (i) $U_b^*=U_b^0=(f_b,C_b)$ is SUC, and
\bqs
1=\sigma^0_b=\sigma^*_b\geq\sigma_b>\tilde \sigma_b; \label{supplylevel1}
\eqs
(ii) if $f_b=s_b<C_b$, another possible solution is that $U_b^*=(C_b,s_b)$ is SOC, $U_b^0=(d_b^0,s_b^0)\neq U_b^*$ with $d_b^0\geq s_b$ and $s_b^0> s_b$, and
\bqs
\sigma^0_b> \sigma^*_b=\sigma_b=\tilde \sigma_b.
\eqs
Therefore, if $U_b^*$ is SUC, $A_*^0=\emptyset$ and $\pi_b^0>0$ or $A_*^0\neq \emptyset$ and $\gamma^0_b(A_*^0)>\theta^0$; if $U_b^*$ is OC, the interior state may not be the same as the stationary state. In addition, for all $b$, $\tilde \sigma_b\leq \sigma_b\leq \sigma_b^* \leq \sigma_b^0$, where all the equality signs hold when $A_*^0=\emptyset$ and $s_b^0=\sum_{a\in A} C_{a\to b} \delta^0_a$, or when $A_*^0\neq \emptyset$ and $\theta^0=\gamma^0_b(A_*^0)$.
\end{lemma}
{\em Proof}. 
\ben
\item When $A_*^0\neq \emptyset$ and $\gamma^0_b(A_*^0)=\theta^0$, then from \refe{entropy:downstream} we have
\bqs
f_b&=&\theta^0 \sum_{a\in A_*^0}  C_{a\to b} +\sum_{\a\in A\setminus A_*^0} C_{\a\to b} \delta^0(\a)
\\ &=& \frac{s_b^0-\sum_{\a\in A\setminus A_*^0} C_{\a\to b} \delta^0(\a)}{\sum_{a\in A_*^0}  C_{a\to b}} \sum_{a\in A_*^0}  C_{a\to b} +\sum_{\a\in A\setminus A_*^0} C_{\a\to b} \delta^0(\a)=s_b^0,
\eqs
If $f_b<s_b$; i.e., $U_b^*=U_b^0=(f_b,C_b)$ is SUC, and from \refe{entropy:downstream} we have $f_b=s_b^0=C_b<s_b$, which is impossible. Thus $f_b=s_b=s_b^0$, and $U_b^*=U_b^0=(C_b,s_b)$ is OC.
\item When $A_*^0=\emptyset$ and $\pi_b^0=0$, then from \refe{entropy:downstream} we have $f_b=s_b^0$. Similarly, we have $f_b=s_b=s_b^0$, and $U_b^*=U_b^0=(C_b,s_b)$ is OC.
\item When $A_*^0\neq \emptyset$ and $\gamma^0_b(A_*^0)>\theta^0$, then from \refe{entropy:downstream} we have $f_b<s_b^0$. If $f_b<s_b$; i.e., $U_b^*=U_b^0=(f_b,C_b)$ is SUC, then $s_b^0=s_b^*=C_b\geq s_b>f_b$, and $1=\sigma^0_b=\sigma^*_b\geq\sigma_b>\tilde \sigma_b$. If $f_b=s_b$; i.e., $U_b^*=(C_b,s_b)$ is OC, and $U_b^0=(d_b^0,s_b^0)\neq U_b^*$ with $d_b^0\geq s_b$ and $s_b^0>s_b$, which is possible only if $s_b<C_b$. In this case,  $\sigma^0_b> \sigma^*_b=\sigma_b=\tilde \sigma_b$.
\item When $A_*^0=\emptyset$ and $\pi_b^0>0$, similarly, there are two types of solutions:  $U_b^*=U_b^0=(f_b,C_b)$ is SUC, then $s_b^0=s_b^*=C_b\geq s_b>f_b$, and $1=\sigma^0_b=\sigma^*_b\geq\sigma_b>\tilde \sigma_b$; $U_b^*=(C_b,s_b)$ is OC, $U_b^0=(d_b^0,s_b^0)\neq U_b^*$ with $d_b^0\geq s_b$ and $s_b^0>s_b$, which is possible only if $s_b<C_b$, and $\sigma^0_b> \sigma^*_b=\sigma_b=\tilde \sigma_b$.
\een
 \eop

\subsection{Solutions of the flux function}

In addition to the interior critical demand function $\theta^0$ in \refe{deftheta}, we also define three other critical demand functions as follows:
\bqs
\theta&=&\min\{\max_{a\in A} \delta_a,\min_b\max_{A_1}\frac{s_b-\sum_{\a\in A\setminus A_1} d_\a\xi_{\a\to b}}{\sum_{a\in A_1} C_a\xi_{a\to b}}\}=\min\{\max_{a\in A} \delta_a, g(\bfdelta,\bfsigma)\},\\
\theta^*&=&\min\{\max_{a\in A} \delta_a^*,\min_b\max_{A_1}\frac{s_b^*-\sum_{\a\in A\setminus A_1} d_\a^*\xi_{\a\to b}}{\sum_{a\in A_1} C_a\xi_{a\to b}}\}=\min\{\max_{a\in A} \delta_a^*, g(\bfdelta^*,\bfsigma^*)\},\\
\tilde \theta&=&\min\{\max_{a\in A} \tilde\delta_a,\min_b\max_{A_1}\frac{s_b-\sum_{\a\in A\setminus A_1} f_\a \xi_{\a\to b}}{\sum_{a\in A_1} C_a\xi_{a\to b}}\}=\min\{\max_{a\in A} \tilde\delta_a, g(\tilde{\bfdelta},\bfsigma)\}.
\eqs
We define the residue supply of link $b$ by $\pi_b=s_b-\sum_{a=1}^m d_a \xi_{a\to b}$.

In the following, we demonstrate that there exists a global flux function, \refe{initialflux}, satisfying \refe{feasiblestationary}, \refe{feasibleinterior}, and \refe{entropy2}. 
We also show that the local flux function in \refe{entropy2} is invariant.

\begin{theorem} \label{thm:final}
For the Riemann problem, we have
\bqn
\theta^0&=&\theta=\theta^*=\tilde \theta,\\
\tilde \delta_a&=&\min\{\delta_a, \theta\}=\min\{\delta_a^0, \theta^0\}=\min\{\delta_a^*, \theta^*\}=\min\{\delta_a, \tilde \theta\}. \label{equivflow}
\eqn
Therefore, the discrete flux function \refe{entropy2} is invariant.
\end{theorem}
{\em Proof}. First, if $\min_b\pi_b\geq 0$, all upstream links are stationary at UC; i.e., $q_a=d_a$ for $a\in A$. Otherwise, from Lemma \ref{lemma:updemand} we have that $\theta^0<\max_a \delta_a\leq \max_a \delta_a^0$, and $A_*=\{a| \delta^0_a>\theta^0\}$ is not empty. From Theorem \ref{thm:cdl}, there exists $b$, such that $\gamma_b^0(A_*)=\theta^0$. Further from Lemma \ref{lemma:downsupply} we have that $f_b=s_b$. However, from \refe{entropy:downstream} we have the following contradiction:
\bqs
f_b&=&\sum_{a\in A_*} \theta^0 C_a+\sum_{\a\in A\setminus A_*} d_\a \xi_{\a\to b}<\sum_{a\in A}d_a \xi_{a\to b}\leq s_b,
\eqs
since $\theta^0<\max_a \delta_a$ and $\pi_b \geq 0$. 
Thus, no upstream links can have SOC stationary states, and from Lemma \ref{lemma:updemand} we have that either $\delta_a^0=\delta_a^*=\delta_a=\tilde \delta_a<\theta^0$ or $\delta_a^0\geq \delta_a^*=\delta_a=\tilde \delta_a=\theta^0$ for any $a\in A$. 
Thus when $A_*=\{a| \delta^0_a>\theta^0\}$ is not empty, then $\theta^0=\max_{a\in A}\delta_a =\max_{a\in A}\delta_a^*=\max_{a\in A}\tilde \delta_a< \max_{a\in A}\delta_a^0$ \footnote{In this case, there exists at least an interior state which is different from the corresponding stationary state.}; otherwise, $\theta^0= \max_{a\in A}\delta_a =\max_{a\in A}\delta_a^*=\max_{a\in A}\tilde \delta_a= \max_{a\in A}\delta_a^0$. \footnote{It is impossible that $\theta^0>\max_{a\in A}\delta_a^0$ from the definition of $\theta^0$.}
In addition, since $s_b^*\geq s_b\geq f_b$, $s_b^*-\sum_{a=1}^m d_a^0 \xi_{a\to b}\geq 0$, and $s_b-\sum_{a=1}^m f_a \xi_{a\to b}$. From Theorem \ref{thm:cdl}, $\theta=\max_{a\in A}\delta_a$, $\theta^*=\max_{a\in A}\delta_a^*$, and $\tilde \theta=\max_{a\in A}\tilde \delta_a$. Therefore, $\theta=\theta^*=\theta^0=\tilde \theta$, and $\tilde \delta_a=\min\{\delta_a, \theta\}=\min\{\delta_a^*, \theta^*\}=\min\{\delta_a^0, \theta^0\}=\min\{\tilde\delta_a, \tilde\theta\}$.

Second, if $\min_b \pi_b <0$, then at least one upstream has a SOC stationary state. Otherwise, $f_a=d_a$ for all $a$, from \refe{entropy:downstream} we have that $f_b=\sum_{a\in A} d_a \xi_{a\to b}>s_b$ when $\pi_b<0$. This contradicts $f_b\leq s_b$ for all $b$. 
From Lemma \ref{lemma:updemand}, $A_*^0=\{a|\theta^0<\delta^0_a\}$ is not empty. 
Thus $\tilde \delta_a=\theta^0$ for $a\in A_*$, and $\delta_\a^0=\delta_\a=\tilde \delta_\a=\delta_\a^*$ for $\a \in A\setminus A_*^0$.
From Theorem \ref{thm:cdl}, we have $\min_{a\in A_*^0} \delta^0_a>\theta^0=\min_{b\in B} \gamma_b^0(A_*^0)\geq \max_{\a\in A\setminus A_*^0} \delta^0_\a$,
where
\bqs
\gamma_b^0(A_*^0)&=&\frac{s_b^0-\sum_{\a\in A\setminus A_*^0} C_{\a\to b} \delta_\a^0}{\sum_{a\in A_*^0} C_{a\to b}}=\frac{C_b\sigma_b^0-\sum_{\a\in A\setminus A_*^0} C_{\a\to b} \delta_\a}{\sum_{a\in A_*^0} C_{a\to b}}.
\eqs
In addition,
\bqs
\gamma_b(A_*^0)&=&\frac{C_b\sigma_b-\sum_{\a\in A\setminus A_*^0} C_{\a\to b} \delta_\a}{\sum_{a\in A_*^0} C_{a\to b}},\\
\gamma_b^*(A_*^0)&=&\frac{s_b^*-\sum_{\a\in A\setminus A_*^0} C_{\a\to b} \delta^*_\a}{\sum_{a\in A_*^0} C_{a\to b}}=\frac{C_b\sigma_b^*-\sum_{\a\in A\setminus A_*^0} C_{\a\to b} \delta_\a}{\sum_{a\in A_*^0} C_{a\to b}},\\
\tilde \gamma_b(A_*^0)&=&\frac{s_b-\sum_{\a\in A\setminus A_*^0} C_{\a\to b} \tilde \delta_\a}{\sum_{a\in A_*^0} C_{a\to b}}=\frac{C_b\sigma_b-\sum_{\a\in A\setminus A_*^0} C_{\a\to b} \delta_\a}{\sum_{a\in A_*^0} C_{a\to b}}.
\eqs
From Lemma \ref{lemma:downsupply}, there exists $b$ such that $\gamma_b^0(A_*^0)=\theta^0$, for which $\sigma_b^0=\sigma_b^*=\sigma_b=\tilde \sigma_b$, and $\theta^0=\gamma_b(A_*^0)=\tilde \gamma_b(A_*^0)=\gamma_b^*(A_*^0)$; for other $b\in B$, we have that $\gamma_b^0(A_*^0)>\theta^0$, and $f_b=\theta^0 \sum_{a\in A_*^0} C_{a\to b}+\sum_{\a\in A\setminus A_*^0} C_{\a\to b} \delta_\a\leq s_b\leq s_b^*$, which leads to $\theta^0\leq \gamma_b(A_*^0)=\tilde \gamma_b(A_*^0)\leq \gamma_b^*(A_*^0)$. Thus we have
\bqs
\min_{b\in B} \gamma_b(A_*^0)=\min_{b\in B} \gamma_b^*(A_*^0)=\min_{b\in B} \tilde \gamma_b(A_*^0)=\min_{b\in B} \gamma_b^0(A_*^0)=\theta^0.
\eqs

If we denote $A_{*}=\{a| \theta^0<\delta_a\}$, which is the set of upstream links with SOC stationary states, then $A_*$ is not empty. From Lemma \ref{lemma:updemand}, we have that $A_*\subseteq A_*^0$ and $\theta^0=\tilde \delta_a=\delta_a=\delta^*_a<\delta^0_a$ for $\a\in A_*^0\setminus A_*$, which is the set of all upstream links with interior states. 
From Lemma \ref{lemma:updemand}, $\gamma_b(A_*)=\gamma_b^*(A_*)=\theta^0$ when $\gamma_b(A_*^0)=\gamma_b^*(A_*^0)=\theta^0$; and $\gamma_b(A_*)=\gamma_b^*(A_*)$ when $\gamma_b(A_*^0)=\gamma_b^*(A_*^0)>\theta^0$. Therefore, $\min_{b\in B} \gamma_b(A_*)=\min_{b\in B} \gamma_b^*(A_*)=\theta^0$. 
From Lemma \ref{lemma:updemand} we have that $\theta^0<\delta_a\leq \delta^*_a$ for $a\in A_*$, and $\theta^0 \geq \delta_\a=\delta^*_\a$ for $\a\in A\setminus A_*$. Thus we have $\min_{a\in A_*} \delta_a^*\geq\min_{a\in A_*} \delta_a>\min_{b\in B} \gamma_b(A_*)=\min_{b\in B} \gamma_b^*(A_*)=\theta^0 \geq \max_{\a\in A\setminus A_*} \delta_\a=\max_{\a\in A\setminus A_*} \delta_\a^*$, which leads to  $\min_{b\in B} \gamma_b(A_*)=\theta=\theta^*=\theta^0$ from Theorem \ref{thm:cdl}.
Since $s_b\geq f_b=\sum_{a=1}^m C_{a\to b} \tilde \delta_a$,  we have from Theorem \ref{thm:cdl} that $\tilde \theta=\max_{a\in A} \tilde \delta_a=\theta^0$ since $A_*^0$ is non-empty.
In addition, when $\a\in A\setminus A_*^0$, $\theta^0\geq \delta_\a^0=\delta_\a=\delta_\a^*=\tilde \theta_\a$; when $a\in A_*^0$, $\theta^0\leq \delta_a\leq \delta_a^*\leq \delta_a^0$. In both cases, \refe{equivflow} is satisfied. 
 \eop

{\em Remark}. Note that Theorem \ref{thm:final} is also true when we switch the min and max operators in definitions of $\theta$, $\theta^0$, $\theta^*$, and $\tilde \theta$. From the proof of Theorem \ref{thm:cdl}, we can see that these values remain the same when $\min_b\pi_b<0$. But if $\min_b\pi_b\geq 0$, the min-max and max-min operators could yield different values, which are not smaller than $\max_a \delta(a)$. A counter example is as follows: $d=[0.5, 0.5]$, $C=[1,1]$, $\pi=\mat{{c} 0.2\\0.2}$, $\xi=\mat{{cc} 0.8 &0.2\\0.2&0.8}$, then $\gamma=\mat{{ccc} 0.75&1.5 &0.7\\1.5&0.75&0.7}$. In this case, $\min_{b}\max_{A_1} \gamma_b(A_1)=1.5$, and $\max_{A_1}\min_{b} \gamma_b(A_1)=0.75$. But both values are greater than 0.5, the maximum demand level. Thus Theorem \ref{thm:final} still holds.

{\em Remark}. From Theorem \ref{thm:final}, we can find a unique flux function $(f_1,\cdots,f_{m+n})=\FF(d_1,\cdots,d_m,s_{m+1},\cdots,s_{m+n})$ defined in \refe{initialflux} and solve the Riemann problem in the following steps:
\ben
\item Calculate $\theta$ from initial conditions in $d_a$, $s_b$, and $\xi_{a\to b}$.
\item Calculate $f_a=\min\{d_a,\theta C_a\}$, and $f_b=\sum_{a\in A} f_a \xi_{a\to b}$.
\item Determine stationary states and interior states\footnote{The interior states may not be unique.} from Corollary \ref{cor:ss-is}.
\een

If we denote $\Theta=g(\bfdelta,\bfsigma)$, then $\Theta=\theta$ when $\min_{b\in B} \pi_b < 0$, and $\Theta\geq \theta$ when $\min_{b\in B} \pi_b \geq 0$. We can see that $f_a=\min\{d_a,\Theta C_a\}$. We define $s_a^+=\Theta C_a$. Then $q_a=\min\{d_a,s_a^+\}$, and $s_a^+$ can be considered as effective downstream supply of upstream link $a$.
We define $\Theta_{-b}$ by $\theta_{-b}=\min_{\b\in B\setminus\{b\}} \max_{A_1\subseteq A} \gamma_\b(A_1)$\footnote{Note that if $b$ is the only downstream link, then $\theta_{-b}$ is set to $1$.}. Obviously $\Theta_{-b}\geq \Theta$. We have the following lemma.
\begin{lemma} \label{lemma:equivalentdemand}
The in-flux of downstream link $b$ can be written as
\bqs
f_b&=&\min\{d_b^-,s_b\},
\eqs
where
\bqs
d_b^-&=&\sum_{a\in A} \min\{d_a,\Theta_{-b} C_a\}\xi_{a\to b}.
\eqs
That is, $d_b^-$ can be considered as effective upstream demand of downstream link $b$. 
\end{lemma}
{\em Proof}. When $\min_{b\in B} \pi_b \geq 0$, from Theorem \ref{thm:cdl} we have that $\Theta\geq \max_{a\in A} \delta _a$, and $\Theta_{-b}\geq \max_{a\in A} \delta _a$. Thus $d_b^-=\sum_{a\in A} d_a \xi_{a\to b}\leq s_b$ since $\pi_b\geq 0$. From Theorem \ref{thm:final} we have $f_b=d_b^-=\min\{d_b^-,s_b\}$.

When $\min_{b\in B} \pi_b <0$, if $\Theta_{-b}=\Theta$, then $f_b=d_b^-=\sum_{a\in A} f_a\xi_{a\to b}\leq s_b$; if $\Theta_{-b}>\Theta$, then from Theorem \ref{thm:cdl} there exists $A_*\neq \emptyset$ such that $\gamma_b(A_*)=\Theta$ and $\min_{a\in A_*} \delta_a >\Theta\geq \max_{\a\in A\setminus A_*} \delta_a$. Therefore, 
\bqs
s_b&=&\Theta \sum_{a\in A_*} C_a \xi_{a\to b} + \sum_{\a\in A\setminus A_*} d_\a \xi_{\a\to b}=\sum_{a\in A} \min\{d_a,\Theta C_a\}\xi_{a\to b}=f_b.
\eqs
In addition, we have
\bqs
f_b&=&\sum_{a\in A} \min\{d_a,\Theta C_a\}\xi_{a\to b}\leq \sum_{a\in A} \min\{d_a,\Theta_{-b} C_a\}\xi_{a\to b}=d_b^-.
\eqs
Therefore, $f_b=\min\{d_b^-,s_b\}$.
 \eop

From Theorem \ref{thm:final} and Lemma \ref{lemma:updemand}, we have the following corollary.
\begin{corollary} \label{cor:upstream}
$U_a^*$ is SOC if and only if $d_a<s_a^+$, and $U_a^*$ is UC if and only if $d_a\geq s_a^+$. In addition, if $d_a>s_a^+$, then $U_a^*=U_a^0=(C_a,q_a)$ is SOC, and
\bqs
1=\delta^0_a=\delta^*_a\geq \delta_a>\tilde \delta_a=\theta;
\eqs
if $d_a<s_a^+$, then $U_a^*=U_a^0=(d_a,C_a)$ is UC, and 
\bqs
\delta^0_a= \delta^*_a=\delta_a=\tilde \delta_a<\theta;
\eqs
if $d_a=s_a^+$, then $U_a^*=(d_a,C_a)$ is UC, $U_a^0=(d_a^0,s_a^0)$ with $d_a^0\geq d_a$ and $s_a^0\geq d_a$, and 
\bqs
\delta^0_a\geq \delta^*_a=\delta_a=\tilde \delta_a=\theta.
\eqs
In the third case, it is possible that $U_a^0\neq U_a^*$, but $U_a^0=U_a^*$ is also a valid solution when $\delta_a^0=\delta_a$.
\end{corollary}

Similarly, from Theorem \ref{thm:final}, Lemma \ref{lemma:downsupply}, and Lemma \ref{lemma:equivalentdemand} we have the following corollary regarding downstream links. 
\begin{corollary} \label{cor:downstream}
$U_b^*$ is SUC if and only if $s_b>d_b^-$, and $U_b^*$ is OC if and only if $s_b\leq d_b^-$. In addition, if $s_b>d_b^-$, then $U_b^*=U_b^0=(q_b,C_b)$ is SUC, and
\bqs
1=\sigma^0_b=\sigma^*_b\geq \sigma_b>\tilde \sigma_b;
\eqs
if $s_b<d_b^-$, then $U_b^*=U_b^0=(s_b,C_b)$ is OC, and 
\bqs
\sigma^0_b= \sigma^*_b=\sigma_b=\tilde \sigma_b;
\eqs
if $s_b=d_b^-$, then $U_b^*=(s_b,C_b)$ is OC, $U_b^0=(s_b^0,s_b^0)$ with $s_b^0\geq s_b$ and $s_b^0\geq s_b$, and 
\bqs
\sigma^0_b\geq \sigma^*_b=\sigma_b=\tilde \sigma_b.
\eqs
In the third case, is is possible that $U_b^0\neq U_b^*$, but $U_b^0=U_b^*$ is also a valid solution when $\sigma^0_b=\sigma_b$.
\end{corollary}

\section{Discussions}
\subsection{Special cases}
 For a linear junction with $m=n=1$, $\xi_{1\to 2}=1$, and $f_1=f_2$. In this case, $\Theta=\frac{s_2}{C_1}$, $\theta=\min\{\delta_1,\frac{s_2}{C_1}\}$, and $f_1=\min\{\theta,\delta_1\}=\min\{\Theta,\delta_1\}$. Thus the effective supply of link 1 is $s_1^+=\Theta C_1=s_2$, and the effective demand of link 2 is $d_2^-=\min\{d_1,\Theta_{-2} C_1\}=d_1$. From Theorem \ref{thm:final} we have $f_1=f_2=\min\{d_1,s_2\}$.
Thus from Corollaries \ref{cor:upstream} and \ref{cor:downstream} we have the following three cases.
\ben
\item When $s_2>d_1$, $d_1<s_1^+$, and $s_2>d_2^-$. Thus $U_1^*=U_1^0=(d_1,C_1)$ is UC, and $U_2^*=U_2^0=(d_1,C_2)$ is SUC. In this case, there is no interior state on either link.
 
\item When $s_2<d_1$, $d_1>s_1^+$, and $s_2<d_s^-$. Thus $U_1^*=U_1^0=(C_1,s_2)$ is SOC, and $U_2^*=U_2^0=(C_2,s_2)$ is OC. In this case, there is no interior state on either link.

\item When $s_2=d_1$, $d_1=s_1^+$, and $s_2=d_s^-$. Thus, $U_1^*=(d_1,C_1)$ is UC, and $U_2^*=(C_2,s_2)$ is OC. In this case, there can be interior states on both links: $U_1^0=(d_1^0,s_1^0)$, where $d_1^0\geq d_1$, and $s_1^0\geq d_1$; $U_2^0=(d_2^0,s_2^0)$, where $d_2^0\geq s_2$, and $s_2^0\geq s_2$. However, since $f_1=\min\{d_1^0,s_1^0\}$, we have either $d_1^0=d_1$ or $s_2^0=s_2$. Thus $U_1^0=U_1^*$ or $U_2^*=U_2^*$; i.e., there can be only one interior state.
\een

For a merging junction with $m>1$ and $n=1$, $\xi_{a\to m+1}=1$, and
\bqs
\Theta&=&\max_{A_1\subset A}\frac{s_{m+1}-\sum_{\a\in A\setminus A_1} d_\a}{\sum_{a\in A_1} C_a}.
\eqs
Thus the effective downstream supply for link $a\in A$ is $s_a^+=\Theta C_a$, and the effective upstream demand for link $m+1$ is $d_{m+1}^-=\sum_{a\in A} d_a$, since $\Theta_{-b}=1$. Then we can solve for stationary and interior states by following Corollaries \ref{cor:upstream} and \ref{cor:downstream}. In particular, when $m=2$, we have
\bqs
\Theta&=&\max\{\frac{s_3-d_2}{C_1},\frac{s_3-d_1}{C_2},\frac{s_3}{C_1+C_2}\},\\
q_1&=&\min\{d_1,\Theta C_1\},\\
q_2&=&\min\{d_2,\Theta C_2\},
\eqs
which is consistent with the fair merge model in \cite{jin2010merge}: 
\bqs
q_1&=&\min\{d_1,\max\{s_3-d_2,\frac{C_1}{C_1+C_2}s_3\}\},\\
q_2&=&\min\{d_2,\max\{s_3-d_1,\frac{C_2}{C_1+C_2}s_3\}\}.
\eqs

 For a diverging junction with $m=1$ and $n>1$, 
\bqs
\Theta&=&\min_{b=2}^{1+n} \frac{s_b}{ C_1 \xi_{1\to b}}.
\eqs
Thus the effective downstream supply for link 1 is $s_1^+=\Theta C_1=\min_{b=2}^{1+n} \frac{s_b}{ \xi_{1\to b}}$, and the effective upstream demand for link $b\in B$ is $d_b^-=\min\{d_1, \min_{\b=2,\b\neq b}^{n+1} \frac{s_\b}{C_1\xi_{1\to \b}}\}\xi_{1\to b}$. Then we can solve for stationary and interior states by following Corollaries \ref{cor:upstream} and \ref{cor:downstream}.
In particular, 
\bqs
q_1&=&\min\{d_1, \Theta C_1\}=\min\{d_1,\min_{b=2}^{1+n} \frac{s_b}{ \xi_{1\to b}}\},
\eqs
which is consistent with the FIFO diverge model \cite{jin2010_diverge}.

For a junction with $m=2$ and $n=2$, we have
\bqs
\Theta&=&\min\{\max\{\frac{s_4-d_2 \xi_{2\to 4}}{C_1 \xi_{1\to 4}},\frac{s_4-d_1\xi_{1\to 4}}{C_2 \xi_{2\to 4}},\frac{s_4}{C_1 \xi_{1\to 4}+C_2 \xi_{2\to 4}}\},\\
&&\max\{\frac{s_4-d_2 \xi_{2\to 4}}{C_1 \xi_{1\to 4}},\frac{s_4-d_1\xi_{1\to 4}}{C_2 \xi_{2\to 4}},\frac{s_4}{C_1 \xi_{1\to 4}+C_2 \xi_{2\to 4}}\} \},\\
f_1&=&\min\{d_1,\Theta C_1\},\\
f_2&=&\min\{d_2,\Theta C_2\}.
\eqs
From Theorem \ref{thm:cdl}, we can see that there are the following scenarios:
\ben
\item Both links 1 and 2 are stationary at UC if and only if $\Theta\geq \max\{\frac{d_1}{C_1},\frac{d_2}{C_2}\}$.
\item Link 1 is stationary at SOC and link 2 at UC if and only if $\frac{d_1}{C_1}>\Theta\geq \frac{d_2}{C_2}$.
\item Both links 1 and 2 are stationary at SOC if and only if $\Theta< \max\{\frac{d_1}{C_1},\frac{d_2}{C_2}\}$.
\een
We can find all stationary and interior states by following Corollaries \ref{cor:upstream} and \ref{cor:downstream}.

\subsection{A simplified framework without interior states}
Within the framework defined in Section \ref{unifiedframework} and a flux function in interior states in \refe{entropy2}, the Riemann solver is well-defined since fluxes can be calculated from Theorem \ref{thm:final}, and stationary and interior states can be determined from Corollaries \ref{cor:upstream} and \ref{cor:downstream}.

From Corollaries \ref{cor:upstream} and \ref{cor:downstream}, interior states can be the same as stationary states in all scenarios, and we can introduce a simplified framework as follows:
\ben
\item Stationary states arise near the junction on all links and satisfy \refe{feasiblestationary}.
\item An entropy condition is defined by the flux function in stationary states: 
\bqs
\theta^*&=&\min\{\max_{a\in A}\delta_a^*,g(\bfdelta^*,\bfsigma^*)\},\\
f_a&=&\min\{d_a^*,\theta^* C_a\}, \forall a\in A\\
f_b&=&\sum_{a\in A} f_a\xi_{a\to b}, \forall b\in B
\eqs
\een
From Theorem \ref{thm:final}, the Riemann problem is uniquely solved and the out-flux of upstream link $a$ is still given by $f_a=\min\{d_a,\theta C_a\}$. Further, Corollaries \ref{cor:upstream} and \ref{cor:downstream} can be used to determine unique stationary states. 

Note that, however, not all CTM flux functions can be used in the simplified framework. For example, with the following flux function proposed in \cite{jin2004network} 
\bsq\label{non-invariant-ff}
\bqn
f_a&=&\min\{d_a^0, \min_{b\in B}\frac{d_a^0}{\sum_{\a\in A} d_\a^0 \xi_{a\to b}}s_b^0 \},\\
f_b&=&\sum_{a\in A} f_a\xi_{a\to b},
\eqn
\esq
it was shown that $f_a=\min\{d_a, \theta C_a\}$, but it is possible that interior states are different from stationary states. In this case, $\FF^0(d_1^0,\cdots,d_m^0,s_{m+1}^0,\cdots,s_{m+n}^0)=\FF(d_1,\cdots,d_m,s_{m+1},\cdots,s_{m+n})=\FF^*(d_1^*,\cdots,d_m^*,s_{m+1}^*,\cdots,s_{m+n}^*)$, but 
\bqs
\FF^0(d_1^0,\cdots,d_m^0,s_{m+1}^0,\cdots,s_{m+n}^0)\neq \FF(d_1^0,\cdots,d_m^0,s_{m+1}^0,\cdots,s_{m+n}^0).
\eqs
That is, $\FF^0(\cdots)\neq \FF(\cdots)$, and it is not invariant.
Therefore, with an invariant flux function, the network kinematic wave model \refe{link-kw} can be defined in the function space of traditional weak solutions. But for a non-invariant flux function, e.g., \refe{non-invariant-ff}, the function space has to be extended to include interior states.

Clearly we have that
\bqs
\RS(\RS(k_1,\cdots,k_{m+n}))=\RS(k_1,\cdots,k_{m+n}).
\eqs
That is, the Riemann solver is consistent in the sense of \cite{garavello2011conservation}.
In addition, all Godunov flux functions can be used as entropy conditions in the simplified framework, in which interior states are the same as stationary states.

\subsection{A stationary junction network}
A junction network is stationary if and only if all initial states are the same as stationary states; i.e., $U_a=U_a^*$ for $a\in A$, and $U_b=U_b^*$ for $b\in B$.
From Corollaries \ref{cor:upstream} and \ref{cor:downstream}, it is possible that $d_a=d_a^*$ ($a\in A$) and $s_b=s_b^*$ ($b\in B$). Therefore, such stationary states always exist.
In addition, if $U_a$ and $U_b$ are stationary, then
\bqs
f_a&=&\min\{d_a,\Theta C_a\}=\min\{d_a,s_a\},\\
f_b&=&\min\{d_b^-,s_b\}=\min\{d_b,s_b\}.
\eqs
From Corollary \ref{cor:upstream} we can see that, if $U_a=(C_a,s_a)$ is SOC, then $\Theta C_a=s_a=q_a$; if $U_a=(d_a,C_a)$ is UC, then $\Theta\geq \delta_a$.
From Corollary \ref{cor:downstream} we can see that, if $U_b=(d_b,C_b)$ is SUC, $d_b^-=q_b<C_b$; if $U_b=(C_b,s_b)$ is OC, then $d_b^-\geq s_b$.

When a junction network is stationary, we can replace any SOC state $(C_a,\theta C_a)$ by the corresponding SUC state $(\theta C_a, C_a)$, and $\theta$ remains the same; but we may not replace a SUC state by the corresponding SOC state. 
Thus for the same fluxes, there can be multiple combinations of stationary states.
Furthermore, if we replace $d_a$ by $q_a=\theta C_a$ for $a\in A_*$, then all upstream links will be stationary at UC, and the critical demand level and all upstream flow-rates remain unchanged. In addition, if all upstream links are stationary at UC, we can replace an OC stationary state on link $b$, $(C_b,q_b)$, by an UC state, $(q_b,C_b)$. 

\section{Conclusion}
The kinematic wave model of network traffic flow, \refe{link-kw}, is thus well-defined by the following rules.
\bi
\item[R1.] The constitutional law: $q_a(x_a,t)=k_a(x_a,t) v_a(x_a,t)$ at any point $x_a$ on link $a$ and time $t$.
\item[R2.] The fundamental diagram of speed-density relation, $v_a(x_a,t)=V_a(k(x_a,t))$, flow-density relation, $q_a(x_a,t)=Q_a(k(x_a,t))\equiv k_a(x_a,t) V_a(k_a(x_a,t))$, demand-density relation, $d_a(x_a,t)=Q_a(\min\{k_{a,c},k_a(x_a,t)\})$, and supply-density relation, $s_a(x_a,t)=Q_a(\max\{k_{a,c},k_a(x_a,t)\})$.
\item[R3.] Traffic conservation: $\pd {k_a} t+\pd {q_a} {x_a}=0$.
\item[R4.] Weak solutions with interior states: the kinematic wave model can have discontinuous weak solutions and interior states at stationary discontinuities.
\item[R5.] Entropy conditions with local flux functions at any junction at $x$: we denote the set of upstream demands  by ${\bf d}(x^-,t)$ and the set of downstream supplies by ${\bf s}(x^+,t)$, then the set of boundary fluxes ${\bf f}(x,t)=\FF^0({\bf d}(x^-,t),{\bf s}(x^+,t))$.
\ei
Thus, if initial conditions in densities and boundary conditions in demands at origins and supplies at destinations are given, \refe{link-kw} can be uniquely solved with the aforementioned five rules. 

This modeling framework makes the entropy condition explicit by using boundary flux functions. In a sense, it is the reverse process of Godunov method, in which flux functions were derived by solving Riemann problems with entropy conditions defined in characteristics or other approaches \cite{lax1972shock}.
It is possible to extend this framework for more complicated situations for multi-class traffic on multi-lane roads or for other types of intersections. The challenges will be related to identifying fundamental diagrams and developing boundary flux functions. Therefore, the Riemann solver can be used to determine whether a flux function is well-defined both mathematically and physically.

In addition, it is possible to extend this framework to study other systems of hyperbolic conservation laws, in which demand and supply functions are well-defined. For example, numerical E-O flux function and other approximate Riemann solvers can be used as entropy conditions when solving the Burgers equation.

In the future, we will be interested in analyzing traffic dynamics in a road network with the help of the new kinematic wave model. Such a Riemann solver and the corresponding kinematic wave model can be used to study many other transportation network problems.

\end {document}